\documentclass[openany, amssymb]{amsart}
\usepackage[margin=1in]{geometry}
\usepackage{mathrsfs,comment}
\usepackage[usenames,dvipsnames]{color}
\usepackage[normalem]{ulem}
\usepackage{url}
\usepackage[all,arc,2cell]{xy}
\UseAllTwocells
\usepackage{enumerate}
\usepackage{amsmath}
\usepackage{tikz-cd}
\usepackage{hyperref}
\usepackage{cleveref}
\usepackage{mathtools}
\usepackage{xcolor}
\usepackage{caption}
\usepackage{subcaption}
\usepackage{graphicx}
\usepackage{bbm}
\newtheorem{thm}{Theorem}[section]
\newtheorem{cor}{Corollary}[section]
\newtheorem{prop}{Proposition}[section]
\newtheorem{lem}{Lemma}[section]

\theoremstyle{definition}
\newtheorem{defn}{Definition}[section]

\newtheorem{exmp}{Example}[section]

\newtheorem{rem}{Remark}[section]

\makeatletter
\newtheorem*{rep@thm}{\rep@title}

\newcommand{\newrepthm}[2]{%
\newenvironment{rep#1}[1]{%
 \def\rep@title{#2 \ref{##1}}%
 \begin{rep@thm}\itshape}%Shoutout https://tex.stackexchange.com/questions/176891/how-to-display-the-text-in-a-new-environment-italic
 {\end{rep@thm}}}

\newcommand{\newrepcor}[2]{%
\newenvironment{rep#1}[1]{%
 \def\rep@title{#2 \ref{##1}}%
 \begin{rep@cor}\itshape}%Shoutout https://tex.stackexchange.com/questions/176891/how-to-display-the-text-in-a-new-environment-italic
 {\end{rep@cor}}}

\makeatother

\makeatother
\newrepthm{thm}{Theorem}
\newrepthm{lemma}{Lemma}
\newrepthm{cor}{Corollary}

\usepackage{tikz,tikz-3dplot}
\tdplotsetmaincoords{80}{45}
\tdplotsetrotatedcoords{-90}{180}{-90}
\usepackage{graphicx} % Required for inserting images

\makeatletter
\let\c@obs=\c@thm
\let\c@cor=\c@thm
\let\c@prop=\c@thm
\let\c@lem=\c@thm
\let\c@prob=\c@thm
\let\c@con=\c@thm
\let\c@conj=\c@thm
\let\c@defn=\c@thm
\let\c@notn=\c@thm
\let\c@notns=\c@thm
\let\c@exmp=\c@thm
\let\c@ax=\c@thm
\let\c@pro=\c@thm
\let\c@ass=\c@thm
\let\c@warn=\c@thm
\let\c@rem=\c@thm
\let\c@sch=\c@thm
\let\c@equation\c@thm
\numberwithin{equation}{section}
\makeatother

\newcommand{\ehr}{\operatorname{ehr}}
\newcommand{\Ehr}{\operatorname{Ehr}}
\newcommand{\Z}{\operatorname{\mathbb{Z}}}
\newcommand{\N}{\operatorname{\mathbb{N}}}
\newcommand{\R}{\operatorname{\mathbb{R}}}
\newcommand{\cone}{\operatorname{cone}}

\newcommand{\juliet}[1]{\textcolor{teal}{#1}}

\usepackage[backend=biber, style=alphabetic, maxnames=99]{biblatex}
\title{On Weighted Ehrhart Series of Order Polytopes}
\author{Evan Huang, Daniel Hwang, David Jin, Juliet Whidden, and Katelyn Zhou}
\date{\today}

\begin{document}

\begin{abstract}
    We derive a combinatorial formula for the integer-point transform of the cone over an order polytope, expressing the transform solely in terms of the linear extensions of the associated poset and their descent statistics. This formula can be specialized to obtain the weighted Ehrhart series of an order polytope with respect to any monomial weight, thus generalizing the work of Stanley in the unweighted case. We then derive an alternative formulation for the weighted Ehrhart series using a stars-and-bars argument. This latter formula yields an algorithm, now implemented in SageMath, for computing the weighted Ehrhart series that is significantly faster than Normaliz's \texttt{WeightedEhrhartSeries} function on posets with large numbers of linear extensions.

    %\textcolor{red}{Get rid of equivalent in equivalent formulation.}
\end{abstract}

\maketitle

\section{Introduction}\label{sec:intro} 

%\textcolor{red}{ALL: KEEP TRACK OF NOTATION!!!!}

%\textcolor{red}{$s$ is degree of $w$, seems to switch to $m$ in Theorem 4.6. }

%\textcolor{red}{$P$ denotes poset, seems to switch to Fubini polynomial in Theorem 2.1.}

%\textcolor{red}{cardinality of poset switches from $n$ to $d$ in Section 4?}

%\textcolor{red}{$S$ is multiple things, including multiset of $[n]$ in Thm 2.1, subset of $[d]$ in Thm 4.6, and $S$ subset in $R^n$ in beginning of Section 2.}

%\textcolor{red}{Replace all $t$ with $m$...}

\subsection{Ehrhart Theory} 
Let \(\mathcal{P} \subset \R^n\) be an $n$-dimensional integral polytope. %\textcolor{red}{Normally, there is no guaranteee that ambient dimension $n$ matches dimension $d$, but in this case, it does. I think we can replace $d$ dimension with $n$, be careful with David's notation.}
In 1967, Ehrhart introduced the function
\[\ehr_{\mathcal{P}}(m) := |m\mathcal{P}\cap \Z^n|\] and proved it is always a polynomial in $m$ of degree $n$, which we refer to as the \textit{Ehrhart polynomial} of $P.$ We then define its associated \textit{Ehrhart series}, which can be written as a rational function of the form 
\[\Ehr_{\mathcal{P}}(z) := \sum_{m\geq 0}\ehr_{\mathcal{P}}(m)z^m =: \frac{h_{\mathcal{P}}^*(z)}{(1 - z)^{n + 1}}.\]  The numerator $h^*_{\mathcal{P}}(z)$ is called the \textit{\(h^*\)-polynomial} of $\mathcal{P},$ a polynomial in $z$ of degree at most $n$. %Hence it is called => which we refer to as
%= \frac{\sum_{i = 0}^n h^*_iz^i}{(1 - z)^{n + 1}} 

We study a generalization of Ehrhart theory in which we assign a weight to each integral point using a polynomial weight function \(w : \R^n\rightarrow \R\) of degree \(s\). The \textit{weighted Ehrhart polynomial}, \textit{weighted Ehrhart series}, and \textit{weighted \(h^*\)-polynomial} are then defined analogously:
\[\ehr_{\mathcal{P}, w}(m) := \sum_{\mathbf{x}\in m\mathcal{P}\cap \Z^n}w(\mathbf{x}),\]
\[\Ehr_{\mathcal{P}, w}(z) := \sum_{m\geq 0}\ehr_{\mathcal{P}, w}(m)z^m =: \frac{h_{\mathcal{P}, w}^*(z)}{(1 - z)^{n + s + 1}}.\]
Similarly, $\ehr_{\mathcal P,w}(m)$ and $h^*_{\mathcal 
P,w}(z)$ are both polynomials, now of degree at most \(n + s\)  \cite{Bajoetal2024}.  Note that using the constant weight function $w=1$ recovers classical Ehrhart theory. 

%= \frac{\sum_{i = 0}^{n + s} h^*_iz^i}{(1 - z)^{n + s+ 1}} 

Much of the interest in Ehrhart theory stems from its ability to encapsulate important geometric information about the associated polytope. For example, the leading coefficient of the weighted Ehrhart polynomial is equal to the integral of the weight over the polytope.  We explore this correspondence through a special class of polytopes called \textit{order polytopes}.

\subsection{Posets and Order Polytopes}

%changed range of tau from [m] to [m-1] so T_{m+1} maps to [m]!

Let \(P\) be a partially-ordered set (poset) on the ground set \([n] := \{1, 2, \dots, n\}\). We will assume that \(P\) is \textit{naturally labeled}, meaning that \(x\leq_P y\Rightarrow x\leq y\). A map \(\tau : [n]\rightarrow \mathbb{N}\) is called \textit{$P$-order-preserving} if it respects the ordering of \(P\), that is,
\[x\leq_P y\quad\Rightarrow\quad \tau(x)\leq\tau(y)\qquad \forall x,y\in [n].\]
We denote the set of $P$-order-preserving maps from \([n]\) to \(\{0,1, \dots, m-1\}\) by \(\mathcal{T}_m(P)\), and just use the terms order-preserving and \(\mathcal{T}_m\) when the poset is implicit.
The \textit{order polynomial} of \(P\), introduced by Stanley \cite{Stanley1986}, is then defined as
\[\Omega_P(m) = |\mathcal T_m|,\qquad m\geq 1.\]
Stanley showed that the order polynomial $\Omega_{P}(m)$ is a polynomial in \(m\) of degree \(n\). Moreover, it is the Ehrhart polynomial of a polytope associated with the corresponding poset $P$ called the \textit{order polytope} $\mathcal{O}(P)$, defined by Stanley in the same paper.

\begin{defn}[\cite{Stanley1986}]
Given a poset $P$ on the set $[n],$ the order polytope  \(\mathcal{O}(P)\subset \R^n\) is defined by the inequalities \begin{align*}
0 \leq x_i \leq 1 & \quad \text{ for all } i\in P ,\\
x_i\leq x_j & \quad \text{ if } i \leq_{P} j.
\end{align*}
\end{defn}

%\juliet{We should pick a notation $\leq_P$ or $\prec_P$.} \textcolor{red}{We vote the former.}

\begin{figure}[hbt!]
    \centering

\begin{tikzpicture}
% Poset 1
  \draw (0,0) node[draw=black,line width=1pt, fill=white!20,circle,inner sep=1pt] {1} -- (0,1) node[draw=black, line width = 1pt, inner sep = 1pt, fill=white!20,circle] {2} -- (0,2) node[draw=black, line width = 1pt, inner sep = 1pt, fill=white!20,circle] {3};

  \draw (0, -0.75) node[rectangle] {$P_1$};

% Poset 2
  \draw (1.5, 0.5) node[draw=black, line width = 1pt, inner sep = 1pt, fill=white!20,circle] {1} -- (2, 1.5) node[draw=black, line width = 1pt, inner sep = 1pt,fill=white!20,circle] {2} -- (2.5, 0.5) node[draw=black, line width = 1pt, inner sep = 1pt,fill=white!20,circle] {3};

  \draw (2, -0.75) node[rectangle] {$P_2$};
  
% Poset 3
  \draw (4, 0) node[draw=black, line width = 1pt, inner sep = 1pt, fill=white!20,circle] {1};
  \draw (4, 1) node[draw=black, line width = 1pt, inner sep = 1pt, fill=white!20,circle] {2};
  \draw (4, 2) node[draw=black, line width = 1pt, inner sep = 1pt, fill=white!20,circle] {3}; 

  \draw (4, -0.75) node[rectangle] {$P_3$};

\end{tikzpicture}
%horizontal space
\hspace{1cm}
%Poset 1
\begin{tikzpicture}%
	[x={(0.707154cm, 0.334138cm)},
	y={(-0.707060cm, 0.334193cm)},
	z={(0.000008cm, 0.881287cm)},
	scale=1.000000,
	back/.style={loosely dotted, thin},
	edge/.style={color=blue!95!black, thick},
	facet/.style={fill=white!95!black,fill opacity=0.800000},
	vertex/.style={inner sep=1pt,circle,draw=green!25!black,fill=green!75!black,thick}]
%
%
%% This TikZ-picture was produced with Sagemath version 10.3
%% with the command: ._tikz_3d_in_3d and parameters:
%% view = [-0.778400000000000, 0.322400000000000, 0.538700000000000]
%% angle = 75.1100000000000
%% scale = 1
%% edge_color = blue!95!black
%% facet_color = blue!95!black
%% opacity = 0.8
%% vertex_color = green
%% axis = False
%%
%% Coordinate of the vertices:
%%
\coordinate (0.00000, 0.00000, 0.00000) at (0.00000, 0.00000, 0.00000);
\coordinate (0.00000, 0.00000, 1.00000) at (0.00000, 0.00000, 1.00000);
\coordinate (0.00000, 1.00000, 1.00000) at (0.00000, 1.00000, 1.00000);
\coordinate (1.00000, 1.00000, 1.00000) at (1.00000, 1.00000, 1.00000);
%%
%%
%% Drawing edges in the back
%%
%%
%%
%% Drawing vertices in the back
%%
%%
%%
%% Drawing the facets
%%
\fill[facet] (0.00000, 1.00000, 1.00000) -- (0.00000, 0.00000, 0.00000) -- (0.00000, 0.00000, 1.00000) -- cycle {};
\fill[facet] (1.00000, 1.00000, 1.00000) -- (0.00000, 0.00000, 0.00000) -- (0.00000, 0.00000, 1.00000) -- cycle {};
\fill[facet] (1.00000, 1.00000, 1.00000) -- (0.00000, 0.00000, 1.00000) -- (0.00000, 1.00000, 1.00000) -- cycle {};
%%
%%
%% Drawing edges in the front
%%
\draw[edge] (0.00000, 0.00000, 0.00000) -- (0.00000, 0.00000, 1.00000);
\draw[edge] (0.00000, 0.00000, 0.00000) -- (0.00000, 1.00000, 1.00000);
\draw[edge] (0.00000, 0.00000, 0.00000) -- (1.00000, 1.00000, 1.00000);
\draw[edge] (0.00000, 0.00000, 1.00000) -- (0.00000, 1.00000, 1.00000);
\draw[edge] (0.00000, 0.00000, 1.00000) -- (1.00000, 1.00000, 1.00000);
\draw[edge] (0.00000, 1.00000, 1.00000) -- (1.00000, 1.00000, 1.00000);
%%
%%
%% Drawing the vertices in the front
%%
\node[vertex] at (0.00000, 0.00000, 0.00000)     {};
\node[vertex] at (0.00000, 0.00000, 1.00000)     {};
\node[vertex] at (0.00000, 1.00000, 1.00000)     {};
\node[vertex] at (1.00000, 1.00000, 1.00000)     {};
\draw (0, 0, -1.00) node[rectangle, line width = 1pt] {$\mathcal{O}(P_1)$};
\end{tikzpicture}
%
%horizontal space
\hspace{0.75cm}
% Poset 2
\begin{tikzpicture}%
	[x={(0.707154cm, 0.334138cm)},
	y={(-0.707060cm, 0.334193cm)},
	z={(0.000008cm, 0.881287cm)},
	scale=1.000000,
	back/.style={loosely dotted, thin},
	edge/.style={color=blue!95!black, thick},
	facet/.style={fill=white!95!black,fill opacity=0.800000},
	vertex/.style={inner sep=1pt,circle,draw=green!25!black,fill=green!75!black,thick}]
%
%
%% This TikZ-picture was produced with Sagemath version 10.3
%% with the command: ._tikz_3d_in_3d and parameters:
%% view = [-0.778400000000000, 0.322400000000000, 0.538700000000000]
%% angle = 75.1100000000000
%% scale = 1
%% edge_color = blue!95!black
%% facet_color = blue!95!black
%% opacity = 0.8
%% vertex_color = green
%% axis = False
%%
%% Coordinate of the vertices:
%%
\coordinate (0.00000, 0.00000, 0.00000) at (0.00000, 0.00000, 0.00000);
\coordinate (0.00000, 0.00000, 1.00000) at (0.00000, 0.00000, 1.00000);
\coordinate (0.00000, 1.00000, 1.00000) at (0.00000, 1.00000, 1.00000);
\coordinate (1.00000, 0.00000, 1.00000) at (1.00000, 0.00000, 1.00000);
\coordinate (1.00000, 1.00000, 1.00000) at (1.00000, 1.00000, 1.00000);
%%
%%
%% Drawing edges in the back
%%
\draw[edge,back] (0.00000, 0.00000, 0.00000) -- (1.00000, 1.00000, 1.00000);
%%
%%
%% Drawing vertices in the back
%%
%%
%%
%% Drawing the facets
%%
\fill[facet] (1.00000, 1.00000, 1.00000) -- (0.00000, 1.00000, 1.00000) -- (0.00000, 0.00000, 1.00000) -- (1.00000, 0.00000, 1.00000) -- cycle {};
\fill[facet] (1.00000, 0.00000, 1.00000) -- (0.00000, 0.00000, 0.00000) -- (0.00000, 0.00000, 1.00000) -- cycle {};
\fill[facet] (0.00000, 1.00000, 1.00000) -- (0.00000, 0.00000, 0.00000) -- (0.00000, 0.00000, 1.00000) -- cycle {};
%%
%%
%% Drawing edges in the front
%%
\draw[edge] (0.00000, 0.00000, 0.00000) -- (0.00000, 0.00000, 1.00000);
\draw[edge] (0.00000, 0.00000, 0.00000) -- (0.00000, 1.00000, 1.00000);
\draw[edge] (0.00000, 0.00000, 0.00000) -- (1.00000, 0.00000, 1.00000);
\draw[edge] (0.00000, 0.00000, 1.00000) -- (0.00000, 1.00000, 1.00000);
\draw[edge] (0.00000, 0.00000, 1.00000) -- (1.00000, 0.00000, 1.00000);
\draw[edge] (0.00000, 1.00000, 1.00000) -- (1.00000, 1.00000, 1.00000);
\draw[edge] (1.00000, 0.00000, 1.00000) -- (1.00000, 1.00000, 1.00000);
%%
%%
%% Drawing the vertices in the front
%%
\node[vertex] at (0.00000, 0.00000, 0.00000)     {};
\node[vertex] at (0.00000, 0.00000, 1.00000)     {};
\node[vertex] at (0.00000, 1.00000, 1.00000)     {};
\node[vertex] at (1.00000, 0.00000, 1.00000)     {};
\node[vertex] at (1.00000, 1.00000, 1.00000)     {};
\draw (0, 0, -1.00) node[rectangle, line width = 1pt] {$\mathcal{O}(P_2)$};
\end{tikzpicture}
%
%horizontal space
\hspace{0.75cm}
% Poset 3
\begin{tikzpicture}%
	[x={(0.707154cm, 0.334138cm)},
	y={(-0.707060cm, 0.334193cm)},
	z={(0.000008cm, 0.881287cm)},
	scale=1.000000,
	back/.style={loosely dotted, thin},
	edge/.style={color=blue!95!black, thick},
	facet/.style={fill=white!95!black,fill opacity=0.800000},
	vertex/.style={inner sep=1pt,circle,draw=green!25!black,fill=green!75!black,thick}]
%
%
%% This TikZ-picture was produced with Sagemath version 10.3
%% with the command: ._tikz_3d_in_3d and parameters:
%% view = [-0.778400000000000, 0.322400000000000, 0.538700000000000]
%% angle = 75.1100000000000
%% scale = 1
%% edge_color = blue!95!black
%% facet_color = blue!95!black
%% opacity = 0.8
%% vertex_color = green
%% axis = False
%%
%% Coordinate of the vertices:
%%
\coordinate (0.00000, 0.00000, 0.00000) at (0.00000, 0.00000, 0.00000);
\coordinate (0.00000, 0.00000, 1.00000) at (0.00000, 0.00000, 1.00000);
\coordinate (0.00000, 1.00000, 0.00000) at (0.00000, 1.00000, 0.00000);
\coordinate (0.00000, 1.00000, 1.00000) at (0.00000, 1.00000, 1.00000);
\coordinate (1.00000, 0.00000, 0.00000) at (1.00000, 0.00000, 0.00000);
\coordinate (1.00000, 0.00000, 1.00000) at (1.00000, 0.00000, 1.00000);
\coordinate (1.00000, 1.00000, 0.00000) at (1.00000, 1.00000, 0.00000);
\coordinate (1.00000, 1.00000, 1.00000) at (1.00000, 1.00000, 1.00000);
%%
%%
%% Drawing edges in the back
%%
\draw[edge,back] (0.00000, 1.00000, 0.00000) -- (1.00000, 1.00000, 0.00000);
\draw[edge,back] (1.00000, 0.00000, 0.00000) -- (1.00000, 1.00000, 0.00000);
\draw[edge,back] (1.00000, 1.00000, 0.00000) -- (1.00000, 1.00000, 1.00000);
%%
%%
%% Drawing vertices in the back
%%
\node[vertex] at (1.00000, 1.00000, 0.00000)     {};
%%
%%
%% Drawing the facets
%%
\fill[facet] (1.00000, 1.00000, 1.00000) -- (0.00000, 1.00000, 1.00000) -- (0.00000, 0.00000, 1.00000) -- (1.00000, 0.00000, 1.00000) -- cycle {};
\fill[facet] (0.00000, 1.00000, 1.00000) -- (0.00000, 0.00000, 1.00000) -- (0.00000, 0.00000, 0.00000) -- (0.00000, 1.00000, 0.00000) -- cycle {};
\fill[facet] (1.00000, 0.00000, 1.00000) -- (0.00000, 0.00000, 1.00000) -- (0.00000, 0.00000, 0.00000) -- (1.00000, 0.00000, 0.00000) -- cycle {};
%%
%%
%% Drawing edges in the front
%%
\draw[edge] (0.00000, 0.00000, 0.00000) -- (0.00000, 0.00000, 1.00000);
\draw[edge] (0.00000, 0.00000, 0.00000) -- (0.00000, 1.00000, 0.00000);
\draw[edge] (0.00000, 0.00000, 0.00000) -- (1.00000, 0.00000, 0.00000);
\draw[edge] (0.00000, 0.00000, 1.00000) -- (0.00000, 1.00000, 1.00000);
\draw[edge] (0.00000, 0.00000, 1.00000) -- (1.00000, 0.00000, 1.00000);
\draw[edge] (0.00000, 1.00000, 0.00000) -- (0.00000, 1.00000, 1.00000);
\draw[edge] (0.00000, 1.00000, 1.00000) -- (1.00000, 1.00000, 1.00000);
\draw[edge] (1.00000, 0.00000, 0.00000) -- (1.00000, 0.00000, 1.00000);
\draw[edge] (1.00000, 0.00000, 1.00000) -- (1.00000, 1.00000, 1.00000);
%%
%%
%% Drawing the vertices in the front
%%
\node[vertex] at (0.00000, 0.00000, 0.00000)     {};
\node[vertex] at (0.00000, 0.00000, 1.00000)     {};
\node[vertex] at (0.00000, 1.00000, 0.00000)     {};
\node[vertex] at (0.00000, 1.00000, 1.00000)     {};
\node[vertex] at (1.00000, 0.00000, 0.00000)     {};
\node[vertex] at (1.00000, 0.00000, 1.00000)     {};
\node[vertex] at (1.00000, 1.00000, 1.00000)     {};
\draw (0, 0, -1.00) node[rectangle, line width = 1pt] {$\mathcal{O}(P_3)$};
\end{tikzpicture}

    \caption{Examples of posets and corresponding order polytopes.}
    \label{fig:order_examples}
\end{figure}
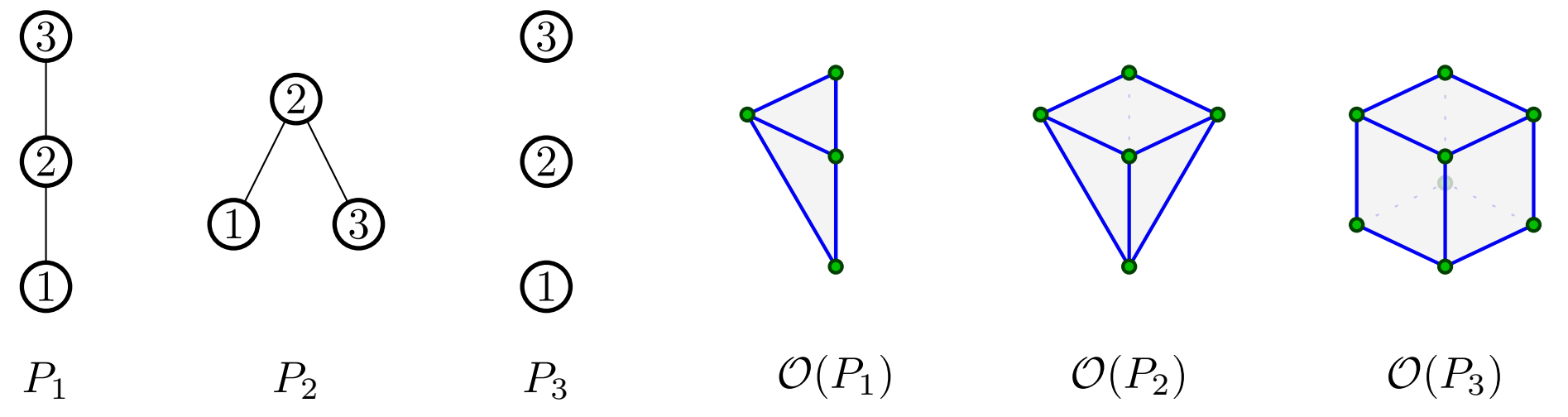

For any poset $P$, $\mathcal{O}(P)$ is an $n$-dimensional integral polytope contained in the unit cube \([0,1]^n\) whose facet-defining hyperplanes are given by the cover relations of \(P\), along with $x_i= 0$ and $x_j= 1$ for all minimal elements $i\in P$ and maximal elements $j\in P.$ Some examples are shown in Figure \ref{fig:order_examples}.

%\juliet{we should redo these images in tikz or Sage without handwritten labels} \textcolor{red}{Do this!}

Moreover, Stanley showed in the same paper that for any poset $P,$ $$\text{ehr}_{\mathcal{O}(P)}(m) = \Omega_{P}(m+1).$$ 
Indeed, if we denote the $m^\text{th}$ dilate of $\mathcal{O}(P)$ by $m\mathcal{O}(P)$, then each integral point $\mathbf{x} \in m\mathcal{O}(P) \cap \mathbb{Z}^n$ respects the ordering on $P$ and satisfies $0 \le x_i \le m$ for all $i \in [n]$. Thus, we can identify $\mathbf x$ with the order-preserving map $\tau \in \mathcal{T}_{m+1}$ given by $\tau(i) = x_i$. 

Stanley then provided a combinatorial formula for the (unweighted) Ehrhart series of \(\mathcal{O}(P)\) expressed in terms of the \textit{linear extensions} of the poset and their \textit{descent numbers}.

%\textcolor{red}{Would it make sense to say that a linear extension $\pi$ of $P$ is an order-preserving bijection from $P$ to itself? And then say $x \leq_P y \implies \pi(x) \leq_P \pi(y)$?}

%\textcolor{red}{As written, I'm pretty sure it should be $i \leq_{\mathbb{N}} j \implies \pi^{-1}(i) \leq_P \pi^{-1}(j)$ below. }

%\textcolor{red}{I'd like to propose to replace first sentence with a linear extension $\pi$ of $P$ is an order preserving bijection from $[n]$ to itself such that $i \leq_P j \implies \pi^{-1}(i) \leq_P \pi^{-1}(j)$}

% , i.e.,  $i \leq_{P} j \implies\pi^{-1}(i) \leq_{\mathbb{N}} \pi^{-1}(j)$

A linear extension $\pi$ of \(P\) is a $P$-order-preserving bijection from \([n]\) to \([n]\), which can also be viewed as a permutation or a word $\pi = (\pi(1), \ldots, \pi(n))$. We denote the set of linear extensions of \(P\) by \(\mathcal{L}(P)\). 
%, and can also be viewed as a permutation or a word $\pi = \pi(1)\dots \pi(n)$
Given a linear extension \(\pi\) of $P$, a \textit{descent} of \(\pi\) is an index \(i\in [n - 1]\) such that \(\pi(i) > \pi(i + 1)\). We use \(d(\pi)\) to denote the number of descents of \(\pi\) and \(\operatorname{Des}(\pi)\) to denote the set of descents of $\pi$. 
For $j\in[n]$, we use $b(\pi, j)$ to denote the number of descents before element $j$ in $\pi$.  More precisely, $$b(\pi, j)=\#\{i \leq \pi^{-1}(j) - 1 :\pi(i)>\pi(i+1)\}.$$
Note that both the order polynomial and the number of linear extensions with a certain number of descents are independent of the natural labeling of the poset.

We now introduce Stanley's combinatorial formula for the (unweighted) Ehrhart series of $\mathcal{O}(P)$:

\begin{repthm}{prop:stanleyorder}[Theorem 4.1, \cite{Stanley1986}; Theorem 3.8.15, \cite{Stanley2012}]
    For any naturally-labeled poset \(P\) on \([n]\),
    \[\Ehr_{\mathcal{O}(P)}(z) = \frac{\sum_{\pi\in\mathcal{L}(P)}z^{d(\pi)}}{(1-z)^{n + 1}}.\]
\end{repthm}

The importance of \Cref{prop:stanleyorder} is manifold. By another result of Stanley in 1980, the \(h^*\)-polynomial of any integral polytope has all non-negative integer coefficients \cite{Stanley1980}. A natural question to ask is whether these coefficients have some combinatorial interpretation in terms of the underlying poset. \Cref{prop:stanleyorder} provides such an interpretation for order polytopes: namely, the coefficient of \(z^k\) is exactly the number of linear extensions of \(P\) with \(k\) descents. Moreover, it provides a way to express the  geometric information of the Ehrhart series purely in terms of poset statistics. For example, it shows that the volume of \(\mathcal{O}(P)\) is given by the number of linear extensions of \(P\) divided by \(n!\).

Motivated by Stanley's unweighted formula, we generalize his work to study the weighted Ehrhart theory of order polytopes.  Note that in the unweighted case, the labeling of the poset does not impact the combinatorial properties of the associated order polytope, as permuting the coordinates does not change the Ehrhart polynomial. However, in the weighted case, the labeling becomes important because of the possible asymmetry of the weight function in the coordinates. 

%\textcolor{red}{Nice work! You go through a lot of information, but the basic goals and background are solid.}

%\textcolor{red}{May I suggest moving the following material on integer-point transform and specialization (up to, not including, Theorem 1.5) to Section 2 and move Outline up here? I say past this point you're probably giving too much information for an introduction - this material more describes how you get Theorem 1.5 and 1.6. Something like ``Our main results are Theorem 1.5 and Theorem 1.6, etc...".}

%\juliet{I agree with this}

Our first main result is an explicit formulation of the weighted Ehrhart series of $\mathcal{O}(P)$: % described in Section \ref{sec:dse formula}:}

\begin{repthm}{thm:weighted_ehrhart_via_differential_operators}
    Let \(P\) be a naturally-labeled poset on \([n]\) and \(w(\mathbf{x}) = \prod_{j\in A}x_j\) be a monomial weight with $A$ a multiset of \([n]\). We define \(\phi : A\rightarrow \{0,1, \dots, n\}\) to be a valid assignment for \(\pi\in \mathcal L(P)\) if for all \(i\), either \(\phi(i) = 0\) or \(\phi(i)\leq\pi^{-1}(i)\). Then
    \[\Ehr_{\mathcal{O}(P), w}(z) = \frac{1}{(1-z)^{n + 1}} \sum_{\pi \in \mathcal{L}(P)} z^{d(\pi)} \sum_{\textup{valid } \phi} \left(\prod_{i\in\phi^{-1}(0)} b(\pi, i)\right) \prod_{k = 1}^n M_{|\phi^{-1}(k)|} \left(\frac{z}{1-z} \right)\]
 where $M_{\ell}(x) = \sum_{j = 0}^{\ell} j! \, S(\ell, j) \, x^j$ is the $\ell^\text{th}$ Fubini polynomial. 
\end{repthm}

%\juliet{What is $P_{|\phi^{-1}(k)|}$ here? } \textcolor{red}{this polynomial is the Fubini polynomial? Find some other notation for it, $F$?} %Define $W_k$ in Section 2?}

%\textcolor{red}{I think it's the name of a polynomial they use later on, it is indeed confusing since $P$ is also the name of our poset. We should change the name of said polynomial. Good catch, rearranged the definitions in Repeat Theorems 2.2 and 2.6.}

The proof of this theorem relies on the following rational expansion of the \textit{integer-point transform} of $\cone(\mathcal{O}(P))$:%, motivated by Stanley's original formula \cite{Stanley2012} in the unweighted case.}

\begin{repthm}{thm:integer-point transform}
Given a naturally-labeled poset \(P\) on \([n]\), 
    \[\sigma_{\cone(\mathcal{O}(P))}(\mathbf{z}) = \frac{1}{1-z_{n + 1}} \sum_{\pi \in \mathcal{L}(P)} z_{n + 1}^{d(\pi)} \left( \prod_{j =1}^n z_j^{b(\pi, j)} \right) \left(\prod_{k=1}^n \frac{1}{1 - W_k(\pi) z_{n + 1}} \right),\]
    where we define $W_k$ by 
    \[W_k(\pi) = 
    \prod_{j=k}^n z_{\pi(j)}
    %\prod_{j \in P  \text{ }\mid \text{ } \pi^{-1}(j) \ge k} z_j
    .\]
\end{repthm}

%\juliet{I forget what we were planning to change the $d(\pi)$ and $b_j(\pi)$ notation to but we should do that.} \textcolor{red}{I think this is resolved?}

%Alternatively, a modified stars-and-bars argument can be used to derive the following equivalent formula. 

Our second main result is an alternative formulation of the weighted Ehrhart series of $\mathcal{O}(P)$ obtained through a modified stars-and-bars argument: % described in Section \ref{sec:stars-and-bars}:

%\juliet{Can we be more consistent with the notation we use for weights? In 2.6 we used $S$ to be the set of indices in the monomial, while in 4.7 we are using $S$ as an arbitrary subset of these indices. Maybe rewrite 2.6 in the notation below and change the domain of $\phi$ to $[d].$} \textcolor{red}{Most natural solution is to change S to T everywhere in Section 4 (Make sure there is no overlap). Kat rewrite weight in Theorem 2.6. Dilation factor is $t$. Degree of weight is $m$. Feel free to override.} \textcolor{red}{Come back to this later.}

%\textcolor{red}{Oct 29: I don't like we are using $s$ for degree weight and $S$ for multiset... didn't we agree to change weight degree to $m$?}

\begin{repthm}{thm:monomial}
Let $P$ be a naturally-labeled poset on the set $[n]$ and let  $w(\mathbf{x})=x_{i_1}x_{i_2}\cdots x_{i_s}$ be any monomial weight.  For each $\pi\in \mathcal{L}(P),$ relabel the $i_1,\dots, i_s$ as $i_1^\pi,\dots, i_s^\pi$ such that $\pi^{-1}(i_1^\pi)\leq\pi^{-1}(i_2^\pi)\leq\cdots\leq\pi^{-1}(i_s^\pi)$. 
For each subset $S\subseteq [s]$, define \[
r(j,S)
=
\#\{\,\ell\in S:\ell<j\,\}.
\]

%\textcolor{red}{Should this be $x_{i_1} \ldots x_{i_s}$ so we are consistent about degree? Also change S to T?}

Then
\[
\begin{aligned}
\operatorname{Ehr}_{\mathcal{O}(P),w}(z)
&=\sum_{\pi\in\mathcal L(P)}
\frac{z^{d(\pi)}}{(1-z)^{n+1}}
\sum_{\substack{S\subseteq [s]}}
\left(\frac{z}{1-z}\right)^{|S|}
\left(
\prod_{j\in S}
\bigl(\pi^{-1}(i_j^\pi)+r(j,S)\bigr)
\right)
\left(
\prod_{j\in [s]\setminus S}
\bigl(b(\pi, i_j^\pi)+r(j,S)\bigr)
\right).
\end{aligned}
\]

%\textcolor{red}{Revisit Theorem 4.5, S notation is inconsistent. I think some $d$ were not changed into $s$. I thought we were changing to $B$?}
\end{repthm}

We can derive from either of the above two results the following combinatorial formula for the integral of any monomial over an order polytope.

%\juliet{Is it okay if I replace the first two sentences of this corollary with ``Using the same notation as in Theorem 4.4, ... "? I don't think we need to restate everything that you can see right above it}

\begin{repcor}{cor: integral combinatorial}
    Using the same notation as above in Theorem \ref{thm:monomial},
    \[\int_{\mathcal{O}(P)}w(\mathbf{x})d\mathbf{x} = \frac{1}{(s + n)!}\sum_{\pi\in\mathcal{L}(P)}\left(\prod_{j\in[n]}\pi^{-1}(i_j^\pi) + (j - 1)\right).\]
\end{repcor}

%\textcolor{red}{David, Juliet, do you have the proof for this?}

Beyond being intrinsically interesting as a result of the purely combinatorial means used to obtain it, the alternative stars-and-bars formulation also yields a much-improved asymptotic runtime complexity when implemented as an algorithm in SageMath \cite{sagemath} to compute the weighted Ehrhart series. Our algorithm consistently outperforms Normaliz' \texttt{WeightedEhrhartSeries} function, the primary freely available computational tool used to compute weighted Ehrhart series, by a factor of at least 10 on posets with between 10 and \(10^5\) elements.

\subsection{Outline}
This paper is organized as follows: In \Cref{sec:dse formula}, we introduce the \textit{integer-point transform} of the cone over a polytope and derive our first formulation for the weighted Ehrhart series of $\mathcal{O}(P)$. In \Cref{sec:ipt formula}, we provide an alternative geometric proof for Theorem \ref{thm:integer-point transform}. In \Cref{sec:stars-and-bars}, we present our second formulation through a modified stars-and-bars argument. In \Cref{sec:algo}, we present information about the performance of our algorithm and compare it to Normaliz's implementation. Finally, we discuss future work in \Cref{sec:future}.

\section{Combinatorial Formula for the Weighted Ehrhart Series}\label{sec:dse formula}

We begin this section by 
expressing the integer-point transform of the cone over the order polytope of a poset in terms of the poset's linear extensions and descent statistics, from which we can recover the weighted Ehrhart series. 

%\textcolor{red}{To ensure consistency with Section 4, I propose to add definitions to integer-point transform and cone. Also note some things defined in Section 4 are defined in Section 2. Definitions are omitted to avoid messing up numbering.}

\begin{defn}
Given a subset \(X\subseteq \mathbb{R}^{n}\), we define the \textit{integer-point transform} of \(X\) as
\[\sigma_X(\mathbf{z}) = \sum_{\mathbf{x}\in X\cap \mathbb{Z}^{n}}\mathbf{z}^\mathbf{x},\] encoding the integer points in \(X\) as a formal sum of Laurent monomials \cite{BeckRobins2015}. 
\end{defn}

\begin{defn}
If a polytope \(\mathcal{P}\subset\R^n\) has vertices \(\mathbf v_1,\dots, \mathbf v_N\), the \textit{cone over} \(\mathcal{P}\) is given by 
\[\operatorname{cone}(\mathcal{P}) = \{\lambda_1( \mathbf v_1, 1) + \dots + \lambda_N(\mathbf v_N, 1)\mid \lambda_i\geq 0\enspace\forall i\in[N]\}\subset \R^{n+1}.\] 
\end{defn}

An example is given below in Figure \ref{fig:cone}.

\begin{figure}[!ht]
\begin{center}

\scalebox{0.75}{
\begin{tikzpicture}[tdplot_main_coords]
  \coordinate (O) at (0,0,0);

  %% make sure to draw everything from back to front
  %\coneback[surface]{-1.5}{2.5}{-15}
  %\coneback[surface]{-3}{2}{-10}
  \draw (0,0,-1) -- (O);

  %\draw[canvas is xy plane at z=4, surface] (-4, 0, 4) -- (-2, 3.4641, 4) -- (2, 3.4641, 4) -- (4, 0, 4)-- (O) --cycle;
  
  %\conefront[surface]{-3}{2}{-10}
  %\conefront[surface]{-1.5}{2.5}{-15}
  %\filldraw[surface] circle (3);
  \draw[->] (-6,0,0) -- (6,0,0) node[right] {$x_1$};
  \draw[->] (0,-6,0) -- (0,6,0) node[right] {$x_2$};
  \draw[->] (O) -- (0,0,5) node[above] {$x_3$};
  
  %\coneback[surface]{1.5}{2.5}{15}
  %\coneback[surface]{3}{2}{10}
  %\conefront[surface]{3}{2}{10}
  %\conefront[surface]{1.5}{2.5}{15}
  %\polyhedralconeback[surface]{2}{2}{3}

  \draw[canvas is xy plane at z=4] (-4, 0, 4) -- (-2, 3.4641, 4) -- (O) --cycle;
  \draw[canvas is xy plane at z=4] (-2, 3.4641, 4) -- (2, 3.4641, 4) -- (O) --cycle;
  \draw[canvas is xy plane at z=4] (2, 3.4641, 4) -- (4, 0, 4)  -- (O) --cycle;
  \draw[canvas is xy plane at z=4] (4, 0, 4)  -- (2, -3.4641, 4) -- (O) --cycle;
  \draw[canvas is xy plane at z=4] (2, -3.4641, 4) -- (-2, -3.4641, 4)  -- (O) --cycle;
  \draw[canvas is xy plane at z=4] (-2, -3.4641, 4) -- (-4, 0, 4)  -- (O) --cycle;

  \draw[canvas is xy plane at z=3] (-3, 0, 3) -- (-1.5, 2.58980, 3) -- (1.5, 2.5980, 3) -- (3, 0, 3) -- (1.5, -2.5980, 3) -- (-1.5, -2.5980, 3) --cycle;
  \draw[canvas is xy plane at z=2] (-2, 0, 2) -- (-1, 1.7320, 2) -- (1, 1.7320, 2) -- (2, 0, 2) -- (1, -1.7320, 2) -- (-1, -1.7320, 2) --cycle;
  \draw[canvas is xy plane at z=1] (-1, 0, 1) -- (-0.5, 0.8660, 1) -- (0.5, 0.8660, 1) -- (1, 0, 1) -- (0.5, -0.8660, 1) -- (-0.5, -0.8660, 1) --cycle;
  
\end{tikzpicture}}
\end{center}
\caption{A visualization of $\cone(\mathcal{P})$ when $\mathcal{P}$ is a hexagon.}

\label{fig:cone}
\end{figure}
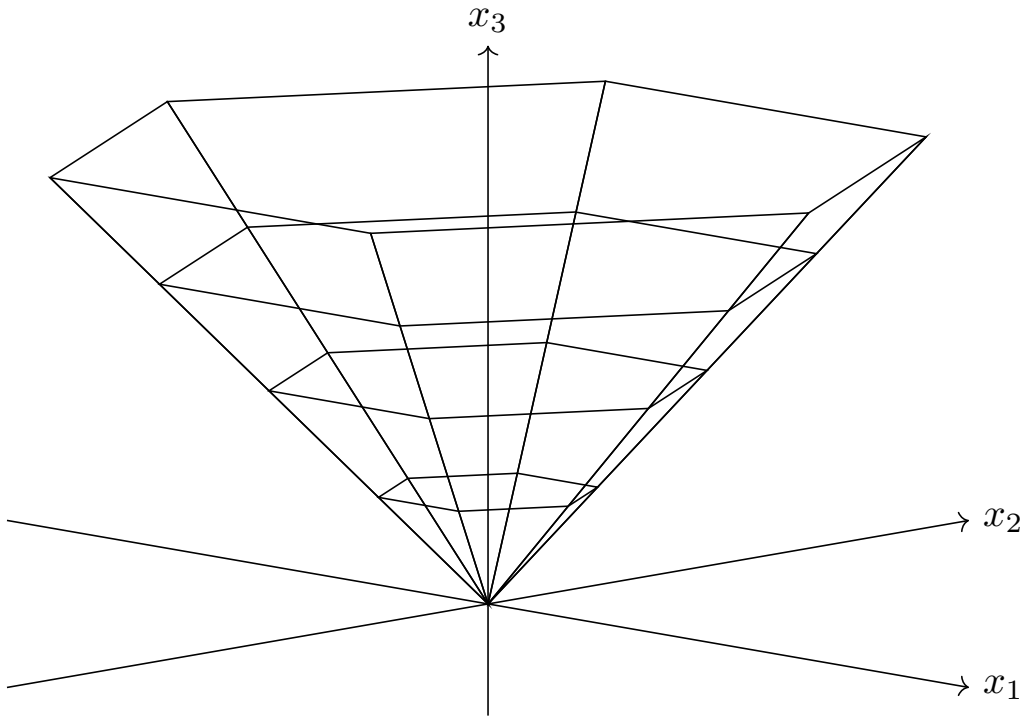

\Cref{thm:integer-point transform} presents our combinatorial closed-form expression for the integer-point transform of the cone over an order polytope. 

Recall from Section \ref{sec:intro} that for a linear extension $\pi\in \mathcal{L}(P)$, $d(\pi)$ is the number of descents of $\pi$ and $b(\pi, j)$ is the number of descents before element $j$ of the poset $P$. 

\begin{thm}\label{thm:integer-point transform}
    Given a naturally-labeled poset \(P\) on \([n]\), 
    \[\sigma_{\cone(\mathcal{O}(P))}(\mathbf{z}) = \frac{1}{1-z_{n + 1}} \sum_{\pi \in \mathcal{L}(P)} z_{n + 1}^{d(\pi)} \left( \prod_{j =1}^n z_j^{b(\pi, j)} \right) \left(\prod_{k=1}^n \frac{1}{1 - W_k(\pi) z_{n + 1}} \right),\]
    where we define $W_k$ by 
    \[W_k(\pi) = \prod_{j=k}^n z_{\pi(j)}
    %\prod_{j \in P  \text{ }\mid \text{ } \pi^{-1}(j) \ge k} z_j
    .\]
\end{thm}

\begin{proof}  

Recall that the integer points of $m\mathcal{O}(P)$ are in bijection with the order-preserving maps $\tau \in \mathcal{T}_{m+1}$. This implies every lattice point of $\cone(\mathcal{O}(P))$ is of the form $(\tau, m)$ for some $\tau \in \mathcal{T}_{m+1}$ and $m \geq 0$ and thus for $\mathbf z = (z_1,\dots, z_{n+1})=(\mathbf z',z_{n+1}),$
\[\sigma_{\operatorname{cone}(\mathcal{O}(P))}(\mathbf{z}) = \sum_{m \ge 0} \sigma_{m\mathcal{O}(P)}(\mathbf{z'}) \, z_{n+1}^m = \sum_{m \ge 0} \sum_{\tau \in \mathcal{T}_{m+1}} \left(\prod_{j=1}^n z_j^{\tau(j)}\right) z_{n+1}^m.\]
    %$$\text{ehr}_{\mathcal{O}(P)}(m) = \Omega_{P}(m+1) = \sum_{\tau \in \mathcal{T}_{m+1}} 1$$

    %We can define the integer-point transform of the $m^{\text{th}}$ dilate $m \mathcal{O}(P)$ as

   % $$\Omega_{P, \mathbf{z}}(m+1) = \sum_{\tau \in \mathcal{T}_{m+1}} \prod_{j=1}^{n} y_{j}^{\tau(j)}$$

   %\textcolor{red}{Observe that all integer lattice points in $\operatorname{cone}(\mathcal{O}(P))$ lie on the hyperplane $z_{n+1} = m$ for some $m \geq 0$. This implies}

    %Lattice points of $\operatorname{cone}(\mathcal{O}(P))$ on the hyperplane $z_{n+1} = m$ are $(\tau, m)$ for $\tau \in \mathcal{T}_{m+1}$, so 
    %\[\sigma_{\operatorname{cone}(\mathcal{O}(P))}(\mathbf{z}) = \sum_{m \ge 0} \sigma_{m\mathcal{O}(P)}(\mathbf{z}) \, z_{n+1}^m = \sum_{m \ge 0} \sum_{\tau \in \mathcal{T}_{m+1}} \left(\prod_{j=1}^n z_j^{\tau(j)}\right) z_{n+1}^m.\]
   % We will expand the right hand side to express it in terms of descent statistics.
    %Note that 
    %\[\sum_{m \ge 0} \Omega_{P, \mathbf{z}}(m+1)z^m = \sum_{m \ge 0} \sum_{\tau \in \mathcal{T}_{m+1}} \left( \prod_{j =1} ^ny_j^{\tau (j)} \right) z^m\]
    %by definition and so what we want to expand out is the right-hand side. 

    %\juliet{I don't think $ \Omega_{P, \mathbf{y}}$ has been defined yet.} \textcolor{red}{Resolved. Switched all $m$ to $t+1$ above, let me know how it reads. Issue with the proof below, should be $\tau \in \mathcal{T}_{m+1}$, not sure if that changes the argument. I think there's an error with the statement below.}
    
    Let $\mathcal T_{\N}$ denote the set of $P$-order-preserving maps $\tau: [n] \rightarrow \mathbb{N}$. Then each $\tau$ appears in the inner summation above exactly when $m \ge \max \tau := \max\limits_{i\in [n]} (\tau(i))$, so we can swap the order of summation to obtain 
    \begin{align*}
        \sum_{m \ge 0} \sum_{\tau \in \mathcal{T}_{m+1}} \left( \prod_{j=1}^n z_j^{\tau (j)} \right) z_{n+1}^m 
        &= \sum_{\tau \in \mathcal T_{\N}} \sum_{m \ge \max \tau} \left( \prod_{j=1}^n z_j^{\tau (j)} \right) z_{n+1}^m\\
        &= \sum_{\tau \in \mathcal T_{\N}} \left( \prod_{j=1}^n z_j^{\tau (j)} \right) \left( \sum_{m \ge \max \tau} z_{n+1}^m\right)\\
        &= \sum_{\tau \in \mathcal T_{\N}} \left( \prod_{j=1}^n z_j^{\tau (j)} \right) \frac{z_{n+1}^{\max \tau}}{1-z_{n+1}}\\
        &= \frac{1}{1-z_{n+1}} \sum_{\tau \in \mathcal T_{\N}} \left( \prod_{j=1}^n z_j^{\tau (j)} \right)z_{n+1}^{\max\tau}.
    \end{align*}

   We now rewrite the summation over $\mathcal T_{\N}$ in terms of linear extensions of $P$. We associate each $\tau \in \mathcal{T}_{\mathbb{N}}$ with the pair $(\pi, \mathbf{t})$, where $\pi$ is the linear extension of $P$ such that $\pi^{-1}(i) < \pi^{-1}(j)$ if and only if either $\tau(i)<\tau(j)$ or $\tau(i)=\tau(j)$ and $i<j$, and   $\mathbf t = (t_1, \dots, t_n) \in \N^n$ is given by $\tau(\pi(1)) = t_1$ and for all $k > 1$, $$\tau(\pi(k)) = \tau(\pi(k-1)) + t_k  + \mathbbm{1}(\pi(k) < \pi(k-1)).$$

   %\textcolor{red}{ Shouldn't this be $\pi(k) <_\mathbb{N} \pi(k-1)$? Explanation follows...}

   %\textcolor{red}{Note the above formula uses $\mathbbm{1}(\pi(k - 1) <_{\mathbb{N}} \pi(k))$ because by our notation above, $k-1 <_P k$ if and only if either $\tau(\pi(k-1)) <_\mathbb{N} \tau(\pi(k))$ or $\tau(\pi(k-1)) = \tau(\pi(k))$ and $\pi(k-1) <_{\mathbb{N}} \pi(k)$. The latter case is impossible whenever $\pi(k) <_{\mathbb{N}} \pi(k-1)$, which means we must enforce $\tau(\pi(k-1)) <_{\mathbb{N}} \tau(\pi(k))$.}
   Note that each descent $\pi(k) < \pi(k-1)$ forces $\tau(\pi(k-1)) < \tau(\pi(k))$, so $t_k\geq 0$ for all $k$.
   Conversely, any pair $(\pi, \mathbf{t})$ defines an order-preserving map $\tau$ whose associated linear extension is $\pi$. This gives a bijection between $\mathcal{T}_{\mathbb{N}}$ and pairs of the form $(\pi, \mathbf{t})$. 
   Using this correspondence, we have 
    \[\tau(j) = \sum_{k = 1}^{\pi^{-1}(j)} t_k + b(\pi, j)\]
    
   %we can recover $\tau$ uniquely from the pair $(\pi,t)$.  Intuitively, we know that $\tau(\pi_k) \le \tau(\pi_{k+1})$ for all $k\in[n-1]$, which is a strict inequality for each descent, and $\mathbf t$ keeps track of the excess distance between adjacent $\tau(\pi_k)$.
    % \textcolor{red}{Note that this argument borrows the idea from Section 4... It almost feels like deja vu...}
    % t represents the difference
    %Since $\pi$ is a linear extension, we know $i \leq_P j \implies \pi(i) \leq_\mathbb{Z} \pi(j)$ 
    
    %Then $d(\pi) = \sum_{k=2}^n \delta_k$ and $b(\pi, i) = \sum_{k = 2}^i \delta_k$. 

    %\textcolor{red}{Skimming heavily through this but looks good. Perhaps give an example of this $\tau \rightarrow (\pi, t)$ mapping?}

    %\textcolor{red}{We are assuming natural labeling and that $\pi_{n} \geq \pi_{i}$ for all $i < n$ by default, right? EDIT: yes. Perhaps make a note of this.}

    %\textcolor{red}{should this be $b(\pi, \pi^{-1}(j))$ everywhere??? Either that, or change your definition of $b(\pi, j)$ to avoid all the rewriting!!! Instead of before index j, it would be before element j. }
    
    and
    \[\max \tau = \tau(\pi(n)) = \sum_{k = 1}^n t_k + d(\pi),\] so we can rewrite the above expression of the integer-point transform as    \begin{align*}\sigma_{\operatorname{cone}(\mathcal{O}(P))}(\mathbf{z})
    & = \frac{1}{1-z_{n+1}} \sum_{\pi\in \mathcal{L}(P)}\sum_{t_1,\dots,t_n\geq 0} \left( \prod_{j=1}^n\left[ z_j^{b(\pi,j)} \prod_{k=1}^{\pi^{-1}(j)} z_j^{t_k}  \right]  \right)\left(z_{n+1}^{d(\pi)} \prod_{k=1}^n z_{n+1}^{t_k}\right)\\
    & = \frac{1}{1-z_{n+1}} \sum_{\pi\in \mathcal{L}(P)} z_{n+1}^{d(\pi)}\left(\prod_{ j =1}^n z_j^{b(\pi, j)}\right)\sum_{t_1,\dots,t_n\geq 0} \prod_{k = 1}^n \Big(W_k(\pi) z_{n+1}\Big)^{t_k}\\
    & = \frac{1}{1-z_{n+1}} \sum_{\pi\in \mathcal{L}(P)} z_{n+1}^{d(\pi)} \left(\prod_{ j =1}^n z_j^{b(\pi, j)}\right) \left( \prod_{k = 1}^n \frac{1}{1 - W_k(\pi)z_{n+1}} \right). \qedhere
    \end{align*}

\end{proof}

\begin{rem}
    \Cref{thm:integer-point transform} can also be derived from Stanley's generating function for $(P, \omega)$-partitions as a special case (Theorem 3.15.5 in \cite{Stanley2012}) for the cone over $\mathcal{O}(P) \times \{1\}$. Our contribution is writing it in the language of linear extension descents, which allows us to differentiate and obtain a formula for the weighted Ehrhart series. 
   % Recall that we can treat lattice points of the cone as order-preserving maps $\tau: \widehat{P} \rightarrow \mathbb{N}$. Then setting $\sigma(a) = \tau(n+2 - a)$ turns them into the order-reversing $(Q, \operatorname{id})$-partitions in Stanley's work. 
    This equivalence can be read in detail in Appendix \ref{appendix:equivalence}.%\textcolor{red}{Make appendix.} %\textcolor{red}{Move section 3 to appendix and update citation here.} %\textcolor{red}{The main idea behind. For more details, please see Appendix.}

\end{rem}

%\begin{rem}
    %Theorem \ref{thm:integer-point transform} as stated gives us reciprocity by a direct application of Stanley's Reciprocity Theorem on the integer cone.  \juliet{can someone clarify this?} 
    
    %\textcolor{red}{Proposed rewrite: Theorem \ref{thm:integer-point transform} as stated gives us reciprocity by a direct application of Stanley's Reciprocity Theorem (Theorem 4.3 of \cite{BeckRobins2015}) on the integer-point transform of $\text{cone}(\mathcal{O}(P))$. From this, we can then recover the weighted Ehrhart series (Proposition 2.1 of \cite{DeLoeraWeighted}).}

    %\textcolor{red}{Honestly, my opinion is to cut this remark.}
%\end{rem}

Our interest in the integer-point transform of the cone over a polytope comes from its ability to recover the weighted Ehrhart series of the polytope with respect to any monomial weight. Given a poset $P$ on $[n]$ and a monomial weight \(w(x_1,\dots, x_n) = \prod_{j\in A}x_j\) for \(A\) a multiset of \([n]\), we have that 
\begin{equation}
\label{eqn:specialization}
\Ehr_{\mathcal{P}, w}(z_{n + 1}) = \left(\prod_{j\in A}\left(z_j\frac{\partial}{\partial z_j}\right)\right)\sigma_{\operatorname{cone}(\mathcal{P})}(z_1,\dots, z_{n + 1})\Big|_{z_1, \dots, z_{n} = 1}. 
\end{equation}
%Got rid of \[ %\], replaced with Equation.

Each $z_j\frac{\partial}{\partial z_j}$ is known as a Cauchy-Euler operator and satisfies $z_j\frac{\partial}{\partial z_j} \mathbf z^\mathbf x= x_j\mathbf z^\mathbf x$, which weighs each point $\mathbf x$ by the monomial weight $w(\mathbf x)=x_j$. 
From the integer-point transform, we obtain the following weighted Ehrhart series formula by applying these Cauchy-Euler operators.

\begin{thm}
\label{thm:weighted_ehrhart_via_differential_operators}
    Let \(P\) be a naturally-labeled poset on \([n]\) and \(w(\mathbf x) = \prod_{j\in A}x_j\) be a monomial weight with $A$ a multiset of \([n]\). We define \(\phi : A\rightarrow \{0,1, \dots, n\}\) to be a valid assignment for \(\pi\in \mathcal L(P)\) if for all \(i\), either \(\phi(i) = 0\) or \(\phi(i)\leq\pi^{-1}(i)\). Then
    \[\Ehr_{\mathcal{O}(P), w}(z) = \frac{1}{(1-z)^{n + 1}} \sum_{\pi \in \mathcal{L}(P)} z^{d(\pi)} \sum_{\textup{valid } \phi} \left(\prod_{i\in\phi^{-1}(0)} b(\pi, i)\right) \prod_{k = 1}^n M_{|\phi^{-1}(k)|} \left(\frac{z}{1-z} \right),\]
 %where $d(\pi)$ is the number of descents in $\pi$, $b(\pi, i)$ is the number of descents in $\pi$ before $i$ and 
 where $M_{\ell}(x) = \sum_{j = 0}^{\ell} j! \, S(\ell, j) \, x^j$ is the $\ell^\text{th}$ Fubini polynomial. 
 %\textcolor{red}{Define some things and then avoid defining them again.}
\end{thm}

\begin{proof}

%We begin with our weighted integer-point transform. Note that in that formula, the `weights' are essentially generating functions of the form $z^{x_1} \dots z^{x_n}$. We can obtain all monomial weights by differentiation -- and the set $A$ that we choose to differentiate is going to be the multiset of our monomial weight. 

For each $i \in \{1, \dots, n\}$, we write
\[D_i = z_i \frac{\partial}{\partial z_i}\]
for the Cauchy-Euler operator, 
and let 
\[D_A = \prod_{i \in A} D_i.\]
Then $D_A\mathbf z^\mathbf x = w(\mathbf x)\mathbf z^\mathbf x,$ so \begin{align*}
D_A\sigma_{\cone(\mathcal{O}(P))}(\mathbf{z})\big|_{z_1,\ldots,z_n=1}
& = \left(\sum_{\mathbf x\in \cone(\mathcal{O}(P))\cap \Z^n} D_A \mathbf z^{\mathbf x}\right)\Bigg|_{z_1,\ldots,z_n=1} = \sum_{\mathbf x\in \cone(\mathcal{O}(P))\cap \Z^n} w(\mathbf x) z_{n+1}^{x_{n+1}},
\end{align*} and the right-hand expression is precisely $\Ehr_{\mathcal{O}(P),w}(z_{n+1}).$  Therefore we can compute the weighted Ehrhart series by applying this differential operator to the formula in Theorem \ref{thm:integer-point transform}, which we restate here:
    \[\sigma_{\cone(\mathcal{O}(P))}(\mathbf{z}) = \frac{1}{1-z_{n + 1}} \sum_{\pi \in \mathcal{L}(P)} z_{n + 1}^{d(\pi)} \left( \prod_{j =1}^n z_j^{b(\pi, j)} \right) \left(\prod_{k=1}^n \frac{1}{1 - W_k(\pi) z_{n + 1}} \right),\quad W_k(\pi) = 
    \prod_{j=k}^n z_{\pi(j)}.
    %\prod_{j \in P  \text{ }\mid \text{ } \pi^{-1}(j) \ge k} z_j.
    \]
Fix a linear extension $\pi\in\mathcal{L}(P)$ and a position $k \in \{1, \dots, n\}$. Let us set 
\[f_k = \frac{1}{1- W_k(\pi)z_{n+1}}\]
and 
\[t_k = \frac{W_k(\pi)z}{1 - W_k(\pi)z_{n+1}}.\]

%\textcolor{red}{Copying Juliet's note, P needs to be defined at the start of Theorem 2.6, also, we need another name for P than Polynomial, here P is reserved for Poset. How about M? or F? IW ILL FIX THIS. I WENT WITH M. }

We note for any $i$ with $\pi^{-1}(i) \ge k$, we have $D_it_k = t_k(1 + t_k)$ and thus, when we apply $D_i$ repeatedly to $f_k$, we can define a polynomial sequence
\[D_{i_1} \dots D_{i_\ell}f_k  = f_k M_\ell(t_k)\]
with $M_0(t) = 1$, $M_1(t) = t$, etc. 
These polynomials satisfy the recurrence
\[M_{\ell+1}(t) = t(1+t)M_\ell'(t) + tM_\ell(t)\]
which gives us the closed-form expression
\[M_\ell(t) = \sum_{j = 0}^d j! S(d,j)t^j\]
where $S(d,j)$ are the Stirling numbers of the second kind. 

%\textcolor{red}{We never formally defined $D_A\sigma_{\text{cone}(\mathcal{O}(P))}(\mathbf{z})$... is what we have fine? } I'm going to do a rewrite that doesn't require this notation. 

%\textcolor{red}{Replace with FPSAC version.}

For a fixed linear extension $\pi\in \mathcal L(P)$, we let $\sigma_{\text{cone}(\mathcal{O}(P))}(\mathbf{z})_{\pi}$ be the integer-point transform of $\text{cone}(\mathcal{O}(P))$ restricted to the points $(\tau, m)$ such that $\pi$ is the linear extension associated to $\tau$ as given in the proof of \Cref{thm:integer-point transform}. By the product rule for derivatives, $D_A\sigma_{\cone(\mathcal{O}(P))}(\mathbf{z})_\pi$ is a sum over maps $\phi: A \rightarrow \{0, 1, \dots, n \}$, where $D_i$ acts on $z_i^{b(\pi, i)}$ when $\phi(i) = 0$ and on $1/(1 - W_k(\pi)z_{n-1})$ if $\phi(i) = k \ge 1$. This holds only for valid assignments as given in the theorem statement, and hence
\[D_A \sigma_{\cone(\mathcal{O}(P))}(\mathbf{z})_\pi = \sigma_{\cone(\mathcal{O}(P))}(\mathbf{z})_\pi \sum_{\text{valid } \phi} \left(\prod_{i \in \phi^{-1}(0)} b(\pi, i) \right) \left( \prod_{k = 1}^n M_{|\phi^{-1}(k)|}(t_k)\right).\]
Taking the sum over $\pi\in \mathcal{L}(P)$ and substituting $z_1, \ldots, z_n = 1$ and  $z_{n+1} = z$ yields the desired result. 
\end{proof}

\section{Geometric Formula for the Integer-Point Transform}\label{sec:ipt formula}

%exercise 5.6 in CRT is the lattice basis
%corollary 6.3.1 are the filters
% homogenizing is ex6.7
%opening faces at descents is 6.2.4.

%\textcolor{red}{Possibly make this a subsection. looks good from heavily skimming}

We offer a geometric proof for the integer-point transform formula in \Cref{thm:integer-point transform}. The proof follows the ideas outlined in Chapter 6 of \cite{Beck_Sanyal_2018}; in particular, the idea of vertices being filter indicators is in Corollary 6.3.1, while the idea of opening faces at descents is in Lemma 6.2.4.%\textcolor{red}{Add numbers, maybe 1 or2 theorems}

\begin{proof}[Alternate Proof to \Cref{thm:integer-point transform}]

Stanley triangulates $\mathcal{O}(P)$ in \cite{Stanley1986} using closed simplices
\[\overline{\Delta_\pi} = \{ x \in \mathbb{R}^n \mid 0 \le x_{\pi(1)} \le \dots \le x_{\pi(n)}\}\quad\text{for }\pi\in\mathcal{L}(P).\]
The corresponding half-open simplices \(\Delta_\pi\) given by converting the nonstrict inequality \(x_{\pi(k)}\leq x_{\pi(k + 1)}\) into a strict inequality for each descent \(k\in[n - 1]\) of $\pi$ yield a disjoint partition of $\mathcal{O}(P)$. Denote by \(C_\pi, \overline{C_\pi}\subseteq\R^{n + 1}\) the cones over \(\Delta_\pi\) and \(\overline{\Delta_\pi}\), respectively. We can compute \(\sigma_{\cone(\mathcal{O}(P))}(\mathbf{z})\) by taking the sum of \(\sigma_{C_\pi}(\mathbf{z})\) over all \(\pi\in\mathcal{L}(P)\). To do this, we first compute the integer-point transform of the \(m^\text{th}\) dilate of \(\overline{\Delta_\pi}\), then sum over all dilates to obtain \(\sigma_{\overline{C_\pi}}(\mathbf{z})\). Then, we deduce \(\sigma_{C_\pi}(\mathbf{z})\) by opening the faces corresponding to descents of \(\pi\), which is equivalent to shifting \(\overline{C_\pi}\) by a subset of its generators. 

Denote the vertex set of \(\overline{\Delta_{\pi}}\) by \(v(\pi) = \{v_0(\pi),v_1(\pi),\dots, v_n(\pi)\}\) and the set of generators for \(\cone(\overline{\Delta_\pi})\) by \(\{u_0(\pi),\dots, u_{n}(\pi)\}\), where $v_0(\pi)=0$ and 
$$v_{k}(\pi) = \sum_{j=k}^n e_{\pi(j)}\quad\text{for }k\in[n],$$
%\juliet{I propose to change this to fix the $P$ vs $[n]$ issue like we did for $W_k$ $$v_{k}(\pi) = \sum_{j=k}^n e_{\pi(j)}\quad\text{for }k\in[n],$$}
$$u_k(\pi) = (v_{k}(\pi), 1)\quad\text{for } 0\le k \le n.$$

Define the $m^\text{th}$ layer of $\overline{\Delta_\pi}$ to consist of all points in the $m^\text{th}$ dilate but not the $(m - 1)^\text{st}$ dilate of $\overline{\Delta_\pi}$. 
Since $\overline{\Delta_\pi}$ is a unimodular simplex,
$\{u_0(\pi), \dots, u_{n}(\pi)\}$ forms a $\mathbb{Z}$-basis of $\mathbb{Z}^{n+1}$.
Therefore each lattice point in the $m^\text{th}$ layer of $\overline{\Delta_\pi}$ can be written uniquely as a nonnegative integer combination of $u_0(\pi), \dots, u_{n}(\pi)$ whose coefficients sum to $m$. Consider a lattice point in $\overline{C}_\pi$ given by $$(\mathbf x,m) = c_0u_0(\pi) + \dots + c_{n}u_{n}(\pi),$$ so  $\mathbf x = c_1v_1(\pi) + \dots + c_n v_n(\pi) \in m\overline{\Delta_\pi}$. 
%Since $v_1(\pi), \dots, v_n(\pi)$ are linearly independent, 
Then $\mathbf x$ is also in  $(m-1)\overline{\Delta_\pi}$ if and only if $c_1 + \dots + c_n \le m-1$, or equivalently, when $c_0 \ge 1$. Thus, the lattice points of the $m^\text{th}$ layer of $\overline{\Delta_\pi}$ are the ones with $c_0 = 0,$ so they correspond bijectively to nonnegative integer combinations $c_1v_1(\pi) + \dots + c_nv_n(\pi)$ with $c_1 + \dots + c_n = m$. 
%There is a bijection between lattice points in the $m^\text{th}$ shell of $\overline{\Delta_\pi}$ and nonnegative integer combinations $\sum_{k = 1}^n c_k v_k(\pi)$ with $\sum c_k = m$. 
Since the integer-point transform of $v_k(\pi)$ is given by $W_k(\pi)$, the integer-point transform of the $m^\text{th}$ layer of $\overline{\Delta_\pi}$ is given by 
%the $m^\text{th}$ coefficient of \(\prod_{k = 1}^n (1 + W_k(\pi) + W_k(\pi)^2 + \dots)\)
%(where the $u_0$-coefficient will be forced to be $0$ in the top shell). 
\[[z^m_{n+1}]\prod_{k = 1}^n\left(\sum_{i\geq 0}W_k(\pi)^iz_{n + 1}^i\right)= [z^m_{n+1}]\prod_{k=1}^n \frac{1}{1 - W_k(\pi) z_{n + 1}},\]
%\textcolor{violet}{Do you think we could cut the first expression so its just the second expression?}
where $[z^m_{n+1}]$ extracts the coefficient of the $z^m_{n+1}$ term.
%the sum of $m$ vectors in $v(\pi)$. 
%Therefore, \[[z^m_{n+1}]\prod_{k=1}^n \frac{1}{1 - W_k(\pi) z_{n + 1}}\] represents the integer-point transform of the $m$th shell of $\overline{\Delta_\pi}$, since each factor $\frac{1}{1- W_k(\pi)z_{n+1}}$ freely sums over the $c_k$ copies of $v_k(\pi)$, each contributing one to the $z_{n+1}$-exponent. 
Since the $m\overline{\Delta_\pi}$ is the disjoint union of its $0^\text{th}$ through $m^\text{th}$ layers of $\overline{\Delta_\pi}$, its integer-point transform is given by 
\begin{align*}
    \sigma_{m\overline{\Delta_\pi}}(\mathbf z) &= \sum_{j=0}^m [z^{j}_{n+1}]\prod_{k=1}^n \frac{1}{1 - W_k(\pi) z_{n + 1}}= [z^m_{n+1}]\frac{1}{1-z_{n + 1}}\prod_{k=1}^n \frac{1}{1 - W_k(\pi) z_{n + 1}}
\end{align*}
where the second equality follows because multiplying by $\frac{1}{1 - z_{n+1}} = \sum_{i \ge 0} z_{n+1}^i$ updates the coefficient of $z_{n+1}^m$ to be the sum of the coefficients of $z_{n+1}^0, \dots, z_{n+1}^m$.
The lattice points of $\overline{C_\pi}$ at height $m$ are the points $(\mathbf{x}, m)$ with $\mathbf{x} \in m\overline{\Delta_\pi} \cap \mathbb{Z}^n$, so we have by summing the integer-point transforms that 
\begin{align*}
    \sigma_{\overline{C_\pi}}(\mathbf{z}) &= \sum_{m \ge 0} \sigma_{m\overline{\Delta_\pi}} (\mathbf{z})\, z_{n+1}^m\\
    &= \sum_{m \ge 0} \left([z^m_{n+1}]\,\frac{1}{1-z_{n + 1}}\prod_{k=1}^n \frac{1}{1 - W_k(\pi) z_{n + 1}}\right)z_{n+1}^m\\
    &= \frac{1}{1-z_{n+1}} \prod_{k=1}^{n} \frac{1}{1-W_k(\pi)z_{n+1}}.
\end{align*}

Now we find \(\sigma_{C_\pi}(\mathbf{z})\) by opening the faces corresponding to descents of \(\pi\). Given a lattice point $(\mathbf{x},m) = \sum_{k=0}^{n} c_k u_k(\pi)$ in \(\overline{C_\pi}\) at height \(m\) where $c_k \in \mathbb{Z}_{\ge 0}$ and $u_k(\pi) = (v_{k}(\pi), 1)$, comparing coordinates yields
\[c_0= m - x_{\pi(n)}, \quad c_1 = x_{\pi(1)}, \quad c_{k+1} = x_{\pi(k+1)} - x_{\pi(k)} \quad (1 \le k \le n-1). \]
A point of $\overline{C_\pi}$ lies in $C_\pi$ if and only if $x_{\pi(k)} < x_{\pi(k+1)}$, or equivalently $c_{k+1} \ge 1$, for every $k \in \operatorname{Des}(\pi)$. Thus,
\[C_\pi \cap \mathbb{Z}^{n+1} = \sum_{k \in \operatorname{Des}(\pi)} u_{k+1} (\pi) + \left( \overline{C_\pi} \cap \mathbb{Z}^{n+1}\right).\]
Since $\mathbf{z}^{u_{k+1}(\pi)} = W_{k+1}(\pi)z_{n+1}$, this gives
%Since $\Delta_{\pi}$ is a simplex, $\overline{C_\pi}$ is unimodular, and we have from our previous result on the integer-point transform of $\sigma_{\overline{C_{\pi}}}$ that%($\sigma$ denoting integer-point transform)
%\[\sigma_{\overline{C_\pi}}(\mathbf{z}) = \prod_{k=1}^{n+1} \frac{1}{1 - W_k(\pi)z_{n+1}}\]

%with $W_{n+1} = 1$ above. Furthermore, we have $x_{\pi(k)} < x_{\pi(k+1)}$ holds if and only if $\lambda_{\pi(k+1)} \ge 1$. So
%\[C_\pi \cap \mathbb{Z}^{n+1} = \sum_{k \in \mathrm{Des}(\pi)} u_{k+1}(\pi) + \left(\overline{C_\pi} \cap \mathbb{Z}^{n+1}\right),\]
%and opening the facet at descent $k$ will shift each lattice point by $u_{k+1}(\pi)$. Since $\mathbf{z}^{u_{k+1}(\pi)} = W_{k+1}(\pi) z_{n+1}$, this yields
\[\sigma_{C_\pi}(\mathbf{z}) = \prod_{k \in \mathrm{Des}(\pi)} W_{k+1}(\pi)z_{n+1} \cdot \sigma_{\overline{C_\pi}}(\mathbf{z}).\]

Each descent contributes one factor of $z_{n+1}$, so this product has degree $d(\pi)$ in $z_{n+1}$. For any given $j \in [n]$, $W_{k+1}(\pi) = \prod_{i = k+1}^n z_{\pi(i)}$ contains \(z_j\) exactly when $k < \pi^{-1}(j)$, so the product in $z_j$ has degree
\[ \#\{k \in \mathrm{Des}(\pi) \mid k < \pi^{-1}(j) \} = b(\pi, j).\]
%\[\#\{k \in \mathrm{Des}(\pi) \mid j \in F_{k+1}(\pi)\} = \#\{k \in \mathrm{Des}(\pi) \mid k < \pi^{-1}(j) \} = b(\pi, j).\]

This gives us
\[\sigma_{C_\pi}(\mathbf{z}) = \frac{1}{1 - z_{n + 1}}\frac{z_{n+1}^{d(\pi)} \prod_{j=1}^n z_j^{b(\pi, j)}}{\prod_{k=1}^{n} \left(1 - W_k(\pi) z_{n+1}\right)}.\]

%Taking the coefficient of $z_{n+1}^m$ is the integer-point transform of $m\Delta_\pi$. Then if we sum over all linear extensions $\pi$, we get the desired result. 

Summing over all linear extensions $\pi$ yields the desired result. 
\end{proof}

\section{Stars-and-Bars Formulation}\label{sec:stars-and-bars}

%\textcolor{red}{Replace all $S$ with $T$.}

With Proposition \ref{exercise}, Stanley provides a formula for the order polynomial of a poset and proves its equivalence to the Ehrhart polynomial of the corresponding order polytope. We present a stars and bars proof of this result in Proposition~\ref{exercise}, and then generalize to the weighted Ehrhart polynomial of order polytopes in Theorem~\ref{thm:monomial}. %We include Lemma~\ref{linear}, an easier subcase of Theorem~\ref{thm:monomial}, in an attempt to bridge the conceptual gap smoothly.

%\begin{defn}
    %Let $P$ be a poset on the set $[n]$. The \emph{order polynomial} $\Omega_{P}(m+1)$ is defined as the number of order-preserving maps $\tau : P \to [0,\dots,m]$. 
%\end{defn}

\begin{defn}
\label{weightdef}
    %Let $P$ be a poset on $n$ elements. An order preserving map $\tau : P \to [0,\dots,m]$ can be represented by $\mathbf{x}=(x_1,\dots,x_n)$, with element $i\in P$ mapped to $x_i\in[0,\dots,m]$, such that $i\leq_P j\implies x_i\leq x_j$. 
    
    For a poset $P$ and polynomial weight $w(\mathbf{x})$, we define the \emph{weighted order polynomial} as $$\Omega_{P,w}(m) = \sum_{\tau \in \mathcal{T}_{m} } w((\tau(1),\dots, \tau(n))),$$ where the summation is taken over order-preserving maps $\tau \in \mathcal{T}_{m}.$ 
\end{defn}

%\textcolor{red}{Daniel: doing quick suggestions in lieu of full edits as of now!}

%\textcolor{red}{This recall definition was already established in Section 2, say something like recall each integral point $x \in m\mathcal{O}(P) \cap \mathbb{Z}^n$ can be identified with the order-preserving map. In particular, this following explanation should be moved into Section 2.}

Recall from Section 1.2 that each integral point $\mathbf{x} \in m\mathcal{O}(P) \cap \mathbb{Z}^n$ can be identified with the order-preserving map $\tau \in \mathcal{T}_{m+1}$ given by $\tau(i) = x_i$. This correspondence gives us the following relationship, with the unweighted case first shown by Stanley (\cite{Stanley1986}).

\begin{lem}
    \label{orderehr}
Let $P$ be a poset on $[n]$ and $w(\mathbf x)$ a polynomial weight function.  Then 
\[\operatorname{ehr}_{\mathcal{O}(P)}(m)=\Omega_P(m+1), \quad \operatorname{ehr}_{\mathcal{O}(P),w}(m)=\Omega_{P,w}(m+1).\]
\end{lem}

\begin{proof}
    As noted above, we identify order-preserving maps $\tau : P \to \{0,\dots,m\}$ with points in $m\mathcal{O}(P)\cap \Z^n,$ so counting lattice points in the $m^\text{th}$ dilate of $\mathcal{O}(P)$ is equivalent to counting maps in $\mathcal{T}_{m+1}$.  Similarly, $\ehr_{\mathcal{O}(P),w}(m)$ weights the coordinates of each lattice point by $w,$ which translates to weighing each map $\tau$ by the weight $w((\tau(1), \ldots, \tau(n)))$, as $\Omega_{P,w}(m+1)$ does.
\end{proof}

\begin{prop}\label{prop:stanleyorder}[Ehrhart Series of an Order Polytope, Theorem 4.1, \cite{Stanley1986}; Theorem 3.8.15, \cite{Stanley2012}]
\label{exercise}
For any naturally-labeled poset $P$ on $[n]$,  
\[\sum_{m\geq 0} \operatorname{ehr}_{\mathcal{O}(P)}(m)z^m=\sum_{m\geq 0} \Omega_P(m+1)z^m=\frac{\sum_{\pi\in \mathcal L(P)} z^{d(\pi)}}{(1-z)^{n+1}}.\]
\end{prop}

\begin{proof}
It suffices to show $$\Omega_P(m+1)=\sum_{\pi\in \mathcal L(P)}\binom{n+m-d(\pi)}{n}.$$

Define $\varphi_m: \mathcal{T}_{m+1}\to \mathcal{L}(P)$ such that $\pi:=\varphi_m(\tau)$ is the linear extension satisfying $\pi^{-1}(i) < \pi^{-1}(j)$ if and only if either $\tau(i)<\tau(j)$ or $\tau(i)=\tau(j)$ and $i<j$.   Then 
$$|\mathcal{T}_{m+1}| = \sum\limits_{\pi\in \mathcal{L}(P)} |\varphi_m^{-1}(\pi)|.$$ 
Let $\pi \in \mathcal{L}(P).$  Each $\tau\in \varphi_m^{-1}(\pi)$ corresponds to a partition of $\{\pi(1),\dots, \pi(n)\}$ into blocks $\tau^{-1}(0),\dots,\tau^{-1}(m)$, each listed in increasing order. Each descent $\pi(k)>\pi(k+1)$ forces $\pi(k)$ and $\pi(k+1)$ to be in different blocks, which we represent by  ``fixed bars.''  We first choose $m-d(\pi)$ ``free bars'' among the $n$ elements of $\pi$ in one of $\binom{n+m-d(\pi)}{n}$ ways.  We then add a ``fixed bar'' at each descent, resulting in a configuration of $n$ stars and $m$ bars which represents some $\tau$ in $\varphi_m^{-1}(\pi)$.
\end{proof}

%\textcolor{red}{I finally understand this proof! Do you mean forces $\pi(i)$ and $\pi(i+1)$ to be in different blocks? You may want to change this index to $k$, multiple usage of index $i$ caused me confusion in parsing your proof. Use $l$ in $S_l = \tau^{-1}(l)$ and write $\pi(k) > \pi(k+1)$ forces $\pi(k)$ and $\pi(k+1)$ to be in different blocks. I propose this change for everywhere you index $\pi$. }

%Let each element of $\pi$ be a star, so we need to separate these $n$ stars into $S_0,\dots,S_m$. In order for each $S_i$ to be in increasing order, whenever there is a descent in $\pi$, we need to create a new $S_i$. Therefore, we have $m-d(\pi)$ bars that we can freely place, then for every descent $\pi(i)\rightarrow\pi(i+1)$, right after the $i$th star we add a fixed bar. Then, we can read off the order preserving map, enclosed between bars $b-1$ and $b$ is stars $x,\dots,y$, then $S_{b}=\{\pi(x),\dots,\pi(y)\}$, which will be in increasing order. We have $n$ stars and $m-d(\pi)$ free bars, so stars and bars gives $\binom{n+m-d(\pi)}{n}$.

%\textcolor{red}{One issue I see is that Sections 1 and 2 do not include many definitions, while Section 4 does. I'm not sure if we want to combine both styles together into one...}
%\textcolor{red}{Say recall here, it's defined in Section 2 already...}
%\textcolor{red}{MOVE MOVE MOVE}
%\begin{defn} For a permutation $\pi\in S_n$ and index $i\in[n]$, let $b(\pi, i)$ be the number of descents before $i$ in $\pi$.  In particular, $b(\pi, i)=\#\{j\in[2,\dots,i]:\pi(j-1)>\pi(j)\}$.  
%\end{defn}

We now give an extension of Stanley's combinatorial formula for $\ehr_{\mathcal{O}(P)}(m)$ in Lemma \ref{prop:stanleyorder} to an analogous combinatorial formula for $\ehr_{\mathcal{O}(P),w}(m)$ for any monomial weight.

\begin{thm}[Monomial Weight Ehrhart Series of Order Polytope]
\label{thm:monomial}
Let $P$ be a naturally-labeled poset on the set $[n]$ and let  $w(\mathbf{x})=x_{i_1}x_{i_2}\cdots x_{i_s}$ be any monomial weight.  For each $\pi\in \mathcal{L}(P),$ relabel the $i_1,\dots, i_s$ as $i_1^\pi,\dots, i_s^\pi$ such that $\pi^{-1}(i_1^\pi)\leq\pi^{-1}(i_2^\pi)\leq\cdots\leq\pi^{-1}(i_s^\pi)$. 
For each subset $S\subseteq [s]$, define \[
r(j,S)
=
\#\{\,\ell\in S:\ell<j\,\}.
\]

Then
\[
\begin{aligned}
\operatorname{Ehr}_{\mathcal{O}(P),w}(z)
&=\sum_{\pi\in\mathcal L(P)}
\frac{z^{d(\pi)}}{(1-z)^{n+1}}
\sum_{\substack{S\subseteq [s]}}
\left(\frac{z}{1-z}\right)^{|S|}
\left(
\prod_{j\in S}
\bigl(\pi^{-1}(i_j^\pi)+r(j,S)\bigr)
\right)
\left(
\prod_{j\in [s]\setminus S}
\bigl(b(\pi,{i_j^\pi})+r(j,S)\bigr)
\right).
\end{aligned}
\]
\end{thm}

\begin{proof}

It suffices to show that $\Omega_{P,w}(m+1)$ is given by 
\begin{equation} \label{weighted order poly formula}
\sum_{\pi\in\mathcal L(P)}
\sum_{\substack{S\subseteq [s]}}
\binom{n+m-d(\pi)}{n+|S|}
\!
\left(
\prod_{j\in S}
\bigl(\pi^{-1}(i_j^\pi)+r(j,S)\bigr)
\right)\!\!
\Bigg(
\prod_{j\in [s]\setminus S}
\bigl(b(\pi,{i_j^\pi})+r(j,S)\bigr)
\Bigg),
\end{equation}
as for any $a \le b$, the coefficient of $z^m$ in $\frac{z^a}{(1-z)^{b+1}}$ is  $\binom{m- a + b}{b}$.

Define $\varphi_m: \mathcal{T}_{m+1}\to \mathcal{L}(P)$ as in Proposition \ref{prop:stanleyorder}. For a fixed $\pi\in \mathcal{L}(P),$ we define certain colored stars and bars configurations such that the number of such configurations is given by the inner sum of (\ref{weighted order poly formula}), and such that each $\tau\in \varphi_m^{-1}(\pi)$ is represented by exactly $w\big((\tau(1),\dots, \tau(n))\big)$ of these configurations.  
Then we would have 
\begin{equation}\label{weighted order poly 2}
\Omega_P(m+1) = \sum_{\pi\in \mathcal{L}(P)}\sum_{\tau\in \varphi_m^{-1}(\pi)}w\big((\tau(1),\dots, \tau(n))\big) = (\ref{weighted order poly formula}), 
\end{equation} as desired. Fix $\pi\in \mathcal{L}(P)$. Take a configuration of $n$ stars, $m-d(\pi)$ ``free bars,'' and $d(\pi)$ ``fixed bars,'' which uniquely represents some $\tau\in \varphi_m^{-1}(\pi)$ as seen in Proposition~\ref{exercise}. Recall that in this context, a ``fixed bar'' is one forced by a descent, otherwise, we have a ``free bar.''

We can weigh each choice of $\tau\in \varphi_m^{-1}(\pi)$ by a linear weight $x_i$ by coloring one bar before the $\pi^{-1}(i)^\text{th}$ star, as there are $\tau(i)$ ways to do so.  
For weight $w(\mathbf{x})=x_{i_1^\pi}\cdots x_{i_s^\pi}$, choose one bar before the $\pi^{-1}(i_c^\pi)^\text{th}$ star to color $c$ for each $c\in[s].$  Assume the colors are distinct and that one bar can hold multiple colors so that these choices are independent.  Then the set $\mathcal{A}_\pi$ of all such colored configurations has size 
$\sum_{\tau\in \varphi_m^{-1}(\pi)}w\big((\tau(1),\dots, \tau(n))\big).$

We now give an alternative construction of $\mathcal{A}_\pi$ counted by (\ref{weighted order poly formula}). 
For a given configuration in $\mathcal A_\pi$, let $C(b)\subseteq[s]$ denote the set of colors present at bar $b$ and let $$S = \big\{\min\big(C(b)\big) : b \text{ is a free bar}\big\}.$$

To count configurations with a given $S=\{c_1,\dots,c_{k}\}\subseteq [s],$ $c_1<\dots<c_{k}$, first choose a configuration of $n+k$ stars and $m-d(\pi)-k$ ``free bars'' in one of $\binom{n+m-d(\pi)}{n+k}$ ways.

We now force minimal colorings at free bars according to $S.$
Choose one of the first $\pi^{-1}(i_{c_k}^\pi)+k-1$ stars, say $s_k$, and turn it into a bar with color $c_k$. Next, choose $s_{k-1}\in [\pi^{-1}(i_{c_{k-1}}^\pi)+k-2]$  and turn the $(s_{k-1})^\text{th}$ star into a bar with color $c_{k-1}$. Continue this process until we have chosen $s_1\in[\pi^{-1}(i_{c_1}^\pi)]$ and turned the $(s_1)^\text{th}$ star into a bar with color $c_{1}$. There are $\prod_{j\in S} (\pi^{-1} (i_j^\pi) + r(j, S))$ ways to do this step. 

For the remaining $j\in [s]\setminus S$, choose one of the $b(\pi,{i_j^\pi})$ fixed bars occurring before the $i_j^{\pi^\text{th}}$ star or one of the $r(j,S)$ already minimally colored free bars occurring before the $i_j^{\pi^\text{th}}$ star to color with color $j$. There are $\prod_{j\in [s]\setminus S}
\bigl(b(\pi,{i_j^{\pi}})+r(j,S)\bigr)$ ways to do this step.

This process generates the same configurations in $\mathcal{A}_\pi$, as all the steps above are reversible. From a configuration of $n$ stars and $m-d(\pi)$ free bars with colors from $S$ distributed on distinct bars such that the only  restriction is that the color $c_j$ must occur on a bar before the $\pi^{-1}(i_{c_j}^\pi)$th star, reversing the transformation gives the corresponding tuple $\{(s_1,s_2,\dots,s_k) \ | \ s_j\in [\pi^{-1}(i_{c_j}^\pi)+j-1],j\in [k]\}$.

Summing over all $S\subset[s]$ and all  $\pi\in \mathcal L(P)$ gives Equality \ref{weighted order poly 2} and completes the proof. \qedhere

\end{proof}

The leading coefficient of the weighted Ehrhart polynomial of an order polytope is equal to the integral of the weight function over the polytope. The only term contributing to the leading coefficient of the weighted Ehrhart polynomial in the formula given in the proof of Theorem~\ref{thm:monomial} is the term when $|S|$ is maximized, i.e. $S=[s]$, which implies the following statement:
\begin{cor}\label{cor: integral combinatorial}
    Let \(P\) be a naturally-labeled poset on \([n]\) and let  $w(\mathbf{x})=x_{i_1}x_{i_2}\cdots x_{i_s}$ be any monomial weight.  For each $\pi\in \mathcal{L}(P),$ relabel the $i_1,\dots, i_s$ as $i_1^\pi,\dots, i_s^\pi$ such that $\pi^{-1}(i_1^\pi)\leq\pi^{-1}(i_2^\pi)\leq\cdots\leq\pi^{-1}(i_s^\pi)$.  Then 
    \[\int_{\mathcal{O}(P)}w(\mathbf{x})d\mathbf{x} = \frac{1}{(s + n)!}\sum_{\pi\in\mathcal{L}(P)}\left(\prod_{j\in[s]}\pi^{-1}(i_j^\pi) + (j - 1)\right).\]
\end{cor}

\section{Algorithmic Implementation}\label{sec:algo}

We have implemented the algorithm for computing the weighted \(h^*\)-polynomial given by \Cref{thm:monomial} in SageMath \cite{sagemath}, available on GitHub \href{https://github.com/Pv2f4n/weighted-order-polytopes}{here}. The motivation for providing an implementation comes from the observation that the asymptotic time complexity of \Cref{thm:monomial} is a significant improvement over that of \Cref{thm:weighted_ehrhart_via_differential_operators}. Whereas \Cref{thm:weighted_ehrhart_via_differential_operators} contains a factor of \(n^s\) due to the iteration over functions \(\phi : A\rightarrow \{0, 1, \dots, n\}\), the former instead substitutes an iteration over subsets of \([s]\), replacing \(n^s\) with \(2^s\). 

\begin{figure}[ht]
    \centering
    \includegraphics[scale=0.45]{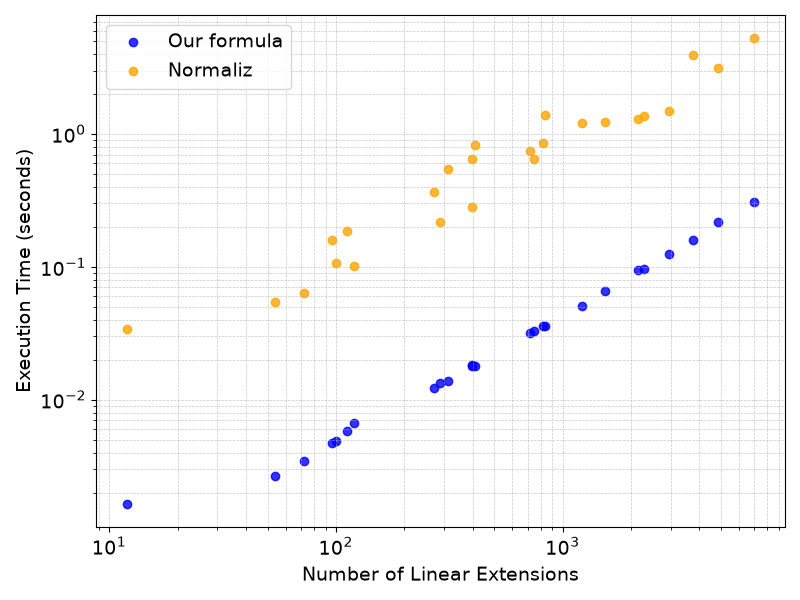}
    \caption{Median-of-three-runs time to compute the weighted \(h^*\)-polynomial with our algorithm (blue) vs \texttt{WeightedEhrhartSeries} (orange) over 20 random posets on 15 elements with \(< 10000\) linear extensions and fixed weight \(w_1(x, y, z) = x^2 - 4xz + 3y\).}
    \label{fig:deg2plot}
\end{figure}

\begin{figure}[ht]
    \centering
    \includegraphics[scale=0.45]{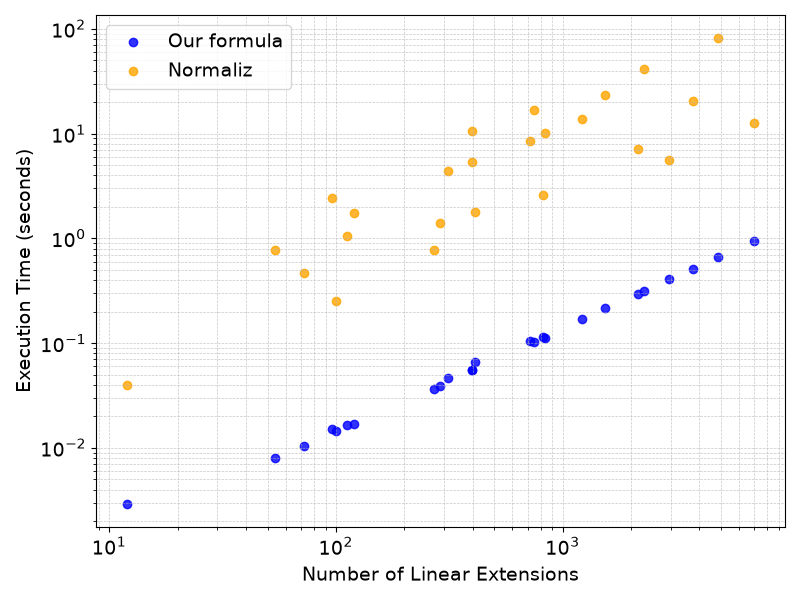}
    \caption{Median-of-three-runs time to compute the weighted \(h^*\)-polynomial with our algorithm (blue) vs \texttt{WeightedEhrhartSeries} (orange) over 20 random posets on 15 elements with \(< 10000\) linear extensions and fixed weight \(w_2(x,y,z) = -3xy^2z - 5y^3z + xy^2 - 2x^2 + 3z + 1\).}
    \label{fig:deg4plot}
\end{figure}

As evidence that this asymptotic improvement is practically useful, we compared our algorithm to the \texttt{WeightedEhrhartSeries} function in Normaliz \cite{normaliz} run through SageMath's PyNormaliz package. We first randomly generated 30 posets on 15 elements where each relation was present with probability 0.5 using the function \texttt{Posets.RandomPoset}, selecting from this set 20 posets with \(<10^5\) linear extensions. We then chose two weights involving three variables, one of degree 2 and the other of degree 4:
\[w_1(x, y, z) = x^2 - 4xz + 3y,\qquad w_2(x,y,z) = -3xy^2z - 5y^3z + xy^2 - 2x^2 + 3z + 1.\]
Next, we compared the runtimes of the two algorithms on each poset with each weight using Python's \texttt{time.perf\_counter} function, taking the median of three runs. We plotted the results on a log-log plot with number of linear extensions versus execution time in seconds using matplotlib \cite{Hunter2007}, as shown in \Cref{fig:deg2plot} and \Cref{fig:deg4plot}. The comparison and plotting functions are also included in the GitHub repository. 

As demonstrated by the plots, our algorithm is consistently between one and two orders of magnitude faster than \texttt{WeightedEhrhartSeries} across both degree 2 and degree 4 weights and varying numbers of linear extensions. In terms of the properties of the poset, our algorithm depends almost exclusively on the number of linear extensions, whereas \texttt{WeightedEhrhartSeries} also shows a significant dependence on the particular structure of the relations in the poset. 

While \texttt{WeightedEhrhartSeries} is a generic algorithm that can be run on any polytope and does not take into account the structure of the underlying poset, our significant practical runtime improvement is nontrivial for two main reasons. First, Normaliz is written in C++ and highly optimized \cite{normaliz}. Second, \texttt{WeightedEhrhartSeries} works by first finding the Stanley decomposition of the polyhedral cone, then integrating over each simplicial piece using Eulerian numbers and Stirling numbers of the second kind \cite{BrunsSoger2015}. For order polytopes, the simplices obtained by Stanley decomposition are exactly those in the canonical unimodular triangulation, and so in a sense \texttt{WeightedEhrhartSeries} already follows the first step of the approach taken by any combinatorial formula. Indeed, implementing the formula in \Cref{thm:weighted_ehrhart_via_differential_operators} directly obtainable from Stanley's \(P\)-partition generating function in SageMath results in a much slower runtime than both \texttt{WeightedEhrhartSeries} and our algorithm even though it does take advantage of poset structure. In this light, we believe that the time improvement from our implementation is both theoretically significant and practically useful.

\section{Future Directions}\label{sec:future}

%\textcolor{red}{David: Good initial draft! You should add one or two sentences at the front about the usefulness of real rootedness. Why are we interested in determining real rootedness? Simple negative real rootedness - log-concave - unimodal - nonnegative Also see stuff below on highlighting results. }

Theorem~\ref{thm:monomial} opens avenues to explore important properties of weighted $h^*$-polynomials within algebraic combinatorics, such as real-rootedness, log-concavity, and unimodality. For a survey of known results and open questions in the unweighted case, see \cite{FerroniHigashitani}. In particular, our theorem gives us a means to represent and efficiently compute weighted $h^*$-polynomials of order polytopes purely in terms of their combinatorial data, allowing us to generate examples and counterexamples and investigate their behavior.

\section{Acknowledgements}\label{sec:acknowledge}

%Written by Daniel.
The second and fourth authors acknowledge support from the U.S. Department of Education Graduate Assistance in Areas of National Need (GAANN) program at Georgia Tech (Award \#P200A240169).

This project began as part of the 2026 Georgia Tech Mathematics REU; it was supported by NSF Grant \#2244427, the Georgia Tech College of Sciences, and the Georgia Tech School of Mathematics. We would like to thank Josephine Yu for her mentorship and helpful feedback. 

\textbf{AI usage declaration}: The results, proofs, and exposition in the final write-up were produced without any AI usage. Gemini 3.8 was used to generate the Python code to convert the raw data from the experiments in \Cref{sec:algo} into matplotlib plots.

% Copied this AI usage declaration from Jochemko and Menon: https://arxiv.org/pdf/2609.17250
%, with the exception of help with the TikZ pictures.

%\textcolor{red}{Acknowledge Georgia Tech School of Math and REU}

%\textcolor{red}{Add AI declaration: We did not use AI in preparing this manuscript.}

{\sloppy\printbibliography}
\newpage
\appendix
\section{Equivalence of Theorem 2.3 to Stanley's P-partitions Generating Function} \label{appendix:equivalence}

   %Let $\widehat{P}$ be the poset with a new maximal element added on top.  
   
   %Recall that we can treat lattice points of the cone as order-preserving maps $\tau: \widehat{P} \rightarrow \mathbb{N}$. %Then setting $\sigma(a) = \tau(n+2 - a)$ turns them into the order-reversing $(Q, \operatorname{id})$-partitions in Stanley's work, where $Q$ is the dual of $\widehat{P}$.
    %\textcolor{red}{Move section 3 to appendix and update citation here.} %\textcolor{red}{The main idea behind. For more details, please see Appendix.}

    %\textcolor{red}{Given this is appendix stuff, I didn't read this part in detail. But what is $F_{P, \omega}$? Is it the integer-point transform somehow? If it is in terms of $(P, \omega)$ partitions, then this should be clearly defined. Also you set $\sigma(a) = \tau(n+2 - a)$, but this should be a definition rather than a setting. In general, if you define anything new, give a brief definition of it.}

    %\textcolor{red}{Define $F_{P, \omega}$ as generating fucntion of $(P, \omega)$-partitions. Briefly define $(P, \omega)$-partitions. Rest of the proof looks good at first glance.}

    %We begin by reintroducing Stanley's notation, all of which can be found in Chapter 15 of \cite{Stanley2012}. 
    
    Let $P$ be a poset naturally labeled by \([p]\).

    \begin{defn}[Definition 3.15.1, \cite{Stanley2012}] A $P$-partition is an order-reversing map $\sigma: P \rightarrow \mathbb{N}$.
    %satisfying the conditions
    % \begin{enumerate}
    %     \item If $s \le t$ in $P$, then $\sigma(s) \ge \sigma(t)$; that is, \(\sigma\) is order-reversing. 
    %     \item If $s < t$ and $\omega(s) > \omega(t)$, then $\sigma (s) > \sigma(t)$. 
    % \end{enumerate}
    The fundamental generating function associated with $P$-partitions is defined by 
    \[F_P(x_1, \dots, x_p) = \sum_{\sigma} x_1^{\sigma(1)} \dots x_p^{\sigma(p)}\]
    for $\sigma$ ranging over all $P$-partitions. 
    \end{defn}

    We record a linear extension \(w\) of $P$ 
    %, regarded as a permutation of the labels $\omega(t)$
    by its word of labels $w=w_1\cdots w_p$.
    %Define $w'_i$ as the index of the element with label $w_i$, i.e.\ $t_{w'_i}=\omega^{-1}(w_i)$.
    The relevant theorem is as follows.

    \begin{thm}[Theorem 3.15.5, \cite{Stanley2012}]\label{thm:stanley p-partitions}
    Let $P$ be a naturally-labeled $p$-element poset. Then
    \[F_P(x_1, \dots, x_p) = \sum_{w \in \mathcal{L}(P)} \frac{\prod_{j \in \operatorname{Des}(w)} x_{w_1}x_{w_2} \dots x_{w_j}}{\prod_{i = 1}^p (1 - x_{w_1}x_{w_2} \dots x_{w_i})}.\]
    %where $\mathcal{L}(P)$ is the set of linear extensions of $P$ and $\operatorname{Des}(w)$ is the descent set of $w$, defined by $D(w) = \{i: w_i > w_{i+1}\}$. 
    \end{thm}

   % For reference, the given theorem states the following: \textit{Let $(P, \omega)$ a labelled $p$-element poset. Then}
    %\[F_{P,\omega}(x_1, \dots, x_p) = \sum_{w \in \mathcal{L}(P, \omega)} \frac{\prod_{j \in D(w)} x_{w_1'}x_{w_2'} \dots x_{w_j'}}{\prod_{i = 1}^p (1 - x_{w_1'}x_{w_2'} \dots x_{w_i'})}.\]

    %  Thus descents are read off the labels, while the \emph{prefix products} $x_{w'_1}\cdots x_{w'_i}$ multiply the variables of the first $i$ elements of the linear extension. Each prefix product gives a factor of the denominator, and those at descents $j\in D(w)$ form the numerator. 

    Recall that lattice points of \(m\mathcal{O}(P)\) are indexed by order-preserving maps \(\tau : P\rightarrow \{0,\dots, m\}\), and so the lattice points of \(\operatorname{cone}(\mathcal{O}(P))\) are indexed by order-preserving maps \(\tau : \widehat{P}\rightarrow\mathbb{N}\) where \(\widehat{P}\) is \(P\) with \(p + 1\) adjoined as the unique maximal element and \(\tau(p + 1)\) represents the height of the lattice point in the cone. Let \(Q\) be the dual of \(\widehat{P}\) under the mapping \(a\mapsto p + 2 - a\). Using \(\tau\) to denote order-preserving maps on \(\widehat{P}\), \(\sigma\) to denote order-reversing maps on \(Q\), and the change of variable \(x_a = z_{p + 2 - a}\), we have 
    \[\sigma_{\operatorname{cone}(\mathcal{O}(P))}(\mathbf{z}) = \sum_\tau z_1^{\tau(1)}\dots z_{p + 1}^{\tau(p + 1)} = \sum_\sigma x_1^{\sigma(1)}\dots x_{p + 1}^{\sigma(p + 1)} = F_Q(x_1,\dots, x_{p + 1}).\]
    By \Cref{thm:stanley p-partitions} applied to \(Q\), 
    \[\sigma_{\operatorname{cone}(\mathcal{O}(P))}(\mathbf{z}) = \sum_{w\in\mathcal{L}(Q)}\frac{\prod_{j\in\operatorname{Des}(w)}x_{w_1}\dots x_{w_j}}{\prod_{i = 1}^{p + 1}(1 - x_{w_1}\dots x_{w_i})} = \sum_{w\in\mathcal{L}(Q)}\frac{\prod_{j\in\operatorname{Des}(w)}z_{p + 2 - w_1}\dots z_{p + 2 - w_j}}{\prod_{i = 1}^{p + 1}(1 - z_{p + 2 - w_1}\dots z_{p + 2 - w_i})}.\]
    Note that there is a bijection 
    \[\ell : \mathcal{L}(\widehat{P})\longrightarrow \mathcal{L}(Q),\qquad \pi_1\dots \pi_{p+1}\longmapsto (p + 2 - \pi_{p+1})\dots (p + 2 - \pi_1),\]
    with $\pi_{p+1}=p+1$.
    Under this bijection, \(w_k = p + 2 - \pi_{p + 2 - k}\), and \(\widehat\pi\) has a descent at \(p + 1 - j\) if and only if \(l(\widehat\pi)\) has a descent at \(j\). So
    \[\sigma_{\operatorname{cone}(\mathcal{O}(P))}(\mathbf{z}) = \sum_{\widehat{\pi}\in\mathcal{L}(\widehat{P})}\frac{\prod_{p + 1 - j\in\operatorname{Des}(\widehat\pi)}z_{\pi_{p + 1}}\dots z_{\pi_{p + 2 - j}}}{\prod_{i = 1}^{p + 1}(1 - z_{\pi_{p + 1}}\dots z_{\pi_{p + 2 - i}})} = \sum_{\pi\in\mathcal{L}(P)}\frac{z^{d(\pi)}_{p + 1}\prod_{j\in\operatorname{Des}(\pi)}z_{\pi_{j + 1}}\dots z_{\pi_p}}{(1 - z_{p + 1})\prod_{k = 1}^p(1 - (z_{\pi_k}\dots z_{\pi_p})z_{p + 1})}\]
    \[= \frac{1}{1 - z_{p + 1}}\sum_{\pi\in\mathcal{L}(P)}z_{p + 1}^{d(\pi)}\frac{\prod_{j = 1}^p z_j^{b(\pi, j)}}{\prod_{k = 1}^p(1 - W_k(\pi)z_{p + 1})}\]
    which is \Cref{thm:integer-point transform}.

\end{document}